\documentclass[12pt]{article}
\usepackage{
amssymb, amsmath}

\title{{\bf Braided tensor categories of admissible modules for affine Lie algebras}}
\author{Thomas Creutzig, Yi-Zhi Huang and Jinwei Yang}
\date{August 9, 2017}
    \begin{document}
    \bibliographystyle{alpha}

\newtheorem{thm}{Theorem}[section]
\newtheorem{defn}[thm]{Definition}
\newtheorem{prop}[thm]{Proposition}
\newtheorem{cor}[thm]{Corollary}
\newtheorem{lemma}[thm]{Lemma}
\newtheorem{rema}[thm]{Remark}
\newtheorem{app}[thm]{Application}
\newtheorem{prob}[thm]{Problem}
\newtheorem{conv}[thm]{Convention}
\newtheorem{conj}[thm]{Conjecture}
\newtheorem{cond}[thm]{Condition}
\newtheorem{criterion}[thm]{Criterion}
    \newtheorem{exam}[thm]{Example}
\newtheorem{assum}[thm]{Assumption}
     \newtheorem{nota}[thm]{Notation}
\newcommand{\halmos}{\rule{1ex}{1.4ex}}
\newcommand{\pfbox}{\hspace*{\fill}\mbox{$\halmos$}}
\newcommand{\nn}{\nonumber \\}

 \newcommand{\res}{\mbox{\rm Res}}
 \newcommand{\ord}{\mbox{\scriptsize \rm ord}}
\renewcommand{\hom}{\mbox{\rm Hom}}
\renewcommand{\exp}{\mbox{\rm exp}}
\newcommand{\edo}{\mbox{\rm End}\ }
\newcommand{\Aut}{\mbox{\rm Aut}\ }
 \newcommand{\pf}{{\it Proof.}\hspace{2ex}}
 \newcommand{\epf}{\hspace*{\fill}\mbox{$\halmos$}}
 \newcommand{\epfv}{\hspace*{\fill}\mbox{$\halmos$}\vspace{1em}}
 \newcommand{\epfe}{\hspace{2em}\halmos}
\newcommand{\nord}{\mbox{\scriptsize ${\circ\atop\circ}$}}
\newcommand{\wt}{\mbox{\rm wt}\ }
\newcommand{\swt}{\mbox{\rm {\scriptsize wt}}\ }
\newcommand{\lwt}{\mbox{\rm wt}^{L}\;}
\newcommand{\rwt}{\mbox{\rm wt}^{R}\;}
\newcommand{\slwt}{\mbox{\rm {\scriptsize wt}}^{L}\,}
\newcommand{\srwt}{\mbox{\rm {\scriptsize wt}}^{R}\,}
\newcommand{\clr}{\mbox{\rm clr}\ }
\newcommand{\tr}{\mbox{\rm Tr}}
\newcommand{\C}{\mathbb{C}}
\newcommand{\Z}{\mathbb{Z}}
\newcommand{\R}{\mathbb{R}}
\newcommand{\Q}{\mathbb{Q}}
\newcommand{\N}{\mathbb{N}}
\newcommand{\CN}{\mathcal{N}}
\newcommand{\F}{\mathcal{F}}
\newcommand{\I}{\mathcal{I}}
\newcommand{\V}{\mathcal{V}}
\newcommand{\one}{\mathbf{1}}
\newcommand{\BY}{\mathbb{Y}}
\newcommand{\ds}{\displaystyle}

        \newcommand{\ba}{\begin{array}}
        \newcommand{\ea}{\end{array}}
        \newcommand{\be}{\begin{equation}}
        \newcommand{\ee}{\end{equation}}
        \newcommand{\bea}{\begin{eqnarray}}
        \newcommand{\eea}{\end{eqnarray}}
         \newcommand{\lbar}{\bigg\vert}
        \newcommand{\p}{\partial}
        \newcommand{\dps}{\displaystyle}
        \newcommand{\bra}{\langle}
        \newcommand{\ket}{\rangle}
\newcommand{\Vir}[1]{\text{Vir}\left(#1\right)}

        \newcommand{\ob}{{\rm ob}\,}
        \renewcommand{\hom}{{\rm Hom}}

\newcommand{\A}{\mathcal{A}}
\newcommand{\Y}{\mathcal{Y}}

\newcommand{\dlt}[3]{#1 ^{-1}\delta \bigg( \frac{#2 #3 }{#1 }\bigg) }

\newcommand{\dlti}[3]{#1 \delta \bigg( \frac{#2 #3 }{#1 ^{-1}}\bigg) }

 \makeatletter
\newlength{\@pxlwd} \newlength{\@rulewd} \newlength{\@pxlht}
\catcode`.=\active \catcode`B=\active \catcode`:=\active
\catcode`|=\active
\def\sprite#1(#2,#3)[#4,#5]{
   \edef\@sprbox{\expandafter\@cdr\string#1\@nil @box}
   \expandafter\newsavebox\csname\@sprbox\endcsname
   \edef#1{\expandafter\usebox\csname\@sprbox\endcsname}
   \expandafter\setbox\csname\@sprbox\endcsname =\hbox\bgroup
   \vbox\bgroup
  \catcode`.=\active\catcode`B=\active\catcode`:=\active\catcode`|=\active
      \@pxlwd=#4 \divide\@pxlwd by #3 \@rulewd=\@pxlwd
      \@pxlht=#5 \divide\@pxlht by #2
      \def .{\hskip \@pxlwd \ignorespaces}
      \def B{\@ifnextchar B{\advance\@rulewd by \@pxlwd}{\vrule
         height \@pxlht width \@rulewd depth 0 pt \@rulewd=\@pxlwd}}
      \def :{\hbox\bgroup\vrule height \@pxlht width 0pt depth
0pt\ignorespaces}
      \def |{\vrule height \@pxlht width 0pt depth 0pt\egroup
         \prevdepth= -1000 pt}
   }
\def\endsprite{\egroup\egroup}
\catcode`.=12 \catcode`B=11 \catcode`:=12 \catcode`|=12\relax
\makeatother

\def\hboxtr{\FormOfHboxtr} 
\sprite{\FormOfHboxtr}(25,25)[0.5 em, 1.2 ex] 

:BBBBBBBBBBBBBBBBBBBBBBBBB | :BB......................B |
:B.B.....................B | :B..B....................B |
:B...B...................B | :B....B..................B |
:B.....B.................B | :B......B................B |
:B.......B...............B | :B........B..............B |
:B.........B.............B | :B..........B............B |
:B...........B...........B | :B............B..........B |
:B.............B.........B | :B..............B........B |
:B...............B.......B | :B................B......B |
:B.................B.....B | :B..................B....B |
:B...................B...B | :B....................B..B |
:B.....................B.B | :B......................BB |
:BBBBBBBBBBBBBBBBBBBBBBBBB |

\endsprite
\def\shboxtr{\FormOfShboxtr} 
\sprite{\FormOfShboxtr}(25,25)[0.3 em, 0.72 ex] 

:BBBBBBBBBBBBBBBBBBBBBBBBB | :BB......................B |
:B.B.....................B | :B..B....................B |
:B...B...................B | :B....B..................B |
:B.....B.................B | :B......B................B |
:B.......B...............B | :B........B..............B |
:B.........B.............B | :B..........B............B |
:B...........B...........B | :B............B..........B |
:B.............B.........B | :B..............B........B |
:B...............B.......B | :B................B......B |
:B.................B.....B | :B..................B....B |
:B...................B...B | :B....................B..B |
:B.....................B.B | :B......................BB |
:BBBBBBBBBBBBBBBBBBBBBBBBB |

\endsprite

\date{}

\maketitle

\begin{abstract}
Using the tensor category theory developed by Lepowsky, Zhang and the second author, we construct a braided tensor category structure with a twist 
on a semisimple category of modules for an affine Lie algebra at an admissible level. We conjecture that this braided tensor category is rigid and thus
is a ribbon category. We also give conjectures on the modularity of this category and on the equivalence with a suitable quantum group tensor category.
In the special case that the affine Lie algebra is $\widehat{\mathfrak{sl}}_2$, we prove the rigidity and modularity conjectures. 
\end{abstract}

\renewcommand{\theequation}{\thesection.\arabic{equation}}
\renewcommand{\thethm}{\thesection.\arabic{thm}}
\setcounter{equation}{0} \setcounter{thm}{0} \date{}
\maketitle

\section{Introduction}
Let $\mathfrak{g}$ be a finite-dimensional simple Lie algebra over $\C$. The affine Lie algebra $\hat{\mathfrak{g}}$ is a central extension of the loop algebra $\mathfrak{g} \otimes \C[t, t^{-1}]$ by a one-dimensional space $\C {\bf k}$. If ${\bf k}$ acts on a $\hat{\mathfrak{g}}$-module by a scalar $\ell \in \C$, we say that the module has level $\ell$. The Bernstein-Gelfand-Gelfand category $\mathcal{O}$ for $\hat{\mathfrak{g}}$ has been extensively studied by both physicists and mathematicians over the years. A central theme is to construct a rigid braided tensor category structure on appropriate subcategories of the category $\mathcal{O}$ at a fixed level.

The works of  Belavin, Polyakov, and Zamolodchikov \cite{BPZ},  Knizhnik and Zamolodchikov \cite{KZ} and 
 Moore and Seiberg \cite{MS} on conformal field theory  led to a conjecture that the category of standard (that is, integrable highest weight) modules at a fixed level $\ell \in \Z_{+}$ has the structure of a modular tensor category. In \cite{KL1}--\cite{KL5}, Kazhdan and Lusztig constructed a natural rigid braided tensor category structure on a certain category of $\hat{\mathfrak{g}}$-modules of level $\ell$ when $\ell + h^{\vee} \notin \Q_{\geq 0}$ (where $h^{\vee}$ is the dual Coexter number of $\mathfrak{g}$ and $\Q_{\geq 0}$ is the set of nonnegative rational numbers), and proved that this rigid braided tensor category is equivalent to a rigid braided tensor category of modules for a quantum group constructed from the same finite-dimensional Lie algebra $\mathfrak{g}$.

Motivated by the work of Kazhdan and Lusztig,  Lepowsky and the second author developed a tensor product theory for the category of modules for a vertex operator algebra in \cite{HL1}--\cite{HL6} and \cite{H1}. Since the category of standard $\hat{\mathfrak{g}}$-modules at a fixed level $\ell \in \Z_{+}$ is the module category for a simple vertex operator algebra $L(\ell,0)$ (\cite{FZ}), they were able to use their theory to construct a braided tensor category structure on the category of finite direct sums of standard modules for $\hat{\mathfrak{g}}$ (\cite{HL7}). Using the method developed in the work of Beilinson, Feigin and Mazur \cite{BFM}, Bakalov and Kirillov \cite{BK} also gave a construction of this braided tensor category structure (but did not give proofs of the rigidity and the nondegeneracy property as is claimed in the book \cite{BK}). In \cite{H6}, using a formula proved and used to derive the Verlinde formula in \cite{H5}, the second author proved the rigidity and nondegeneracy property of this braided tensor category and thus proved the conjecture of Moore-Seiiberg mentioned above that this category has the structure of a modular tensor category.  

Assuming the existence of the rigid braided tensor category structure on the category of standard $\hat{\mathfrak{g}}$-modules at a positive integral level $\ell \in \Z_{+}$  and excluding some cases (in particular, the $\mathfrak{g}=E_{8}$ and $\ell=2$ case), Finkelberg  \cite{F1, F2, F3} proved that this rigid braided tensor categroy
is equivalent to a semisimple subquotient of a rigid braided tensor category of modules for a quantum group constructed from $\mathfrak{g}$, using the equivalence constructed by Kazhdan-Lusztig \cite{KL5} and the Verlinde formula proved by Faltings \cite{Fa}, Teleman \cite{T} and the second author \cite{H5}. This work can also be interpreted as giving another construction of the rigid braided tensor category structure on the category of standard $\hat{\mathfrak{g}}$-modules at a fixed level $\ell \in \Z_{+}$. In \cite{Mc}, motivated by finding a direct proof of Finkelberg's equivalence of tensor categories, McRae gave a direct proof that this tensor category is equivalent to a certain tensor category of finite-dimensional $\mathfrak{g}$-modules.

Using a more general logarithmic tensor category theory for suitable generalized modules for a vertex operator algebra developed by Lepowsky, Zhang and the second author in \cite{HLZ0}--\cite{HLZ8}, Zhang gave in \cite{Zh} (with a mistake corrected in \cite{H4}) a vertex-operator-algebraic construction of the braided tensor categories constructed in \cite{KL1}--\cite{KL5} in the case that $\ell + h^{\vee} \notin \Q_{\geq 0}$.

In the case that $\ell + h^{\vee} \in \Q_{\geq 0}\setminus \Z_{+}$, the category of finitely-generated generalized modules for $L(\ell,0)$ is conjectured to have a natural structure of a braided tensor category. In this paper, we study the case when $\ell$ is an admissible number introduced by Kac and Wakimoto in \cite{KW1} in order to describe and classify admissible $\hat{\mathfrak{g}}$-modules, whose modified characters are holomorphic modular forms with respect to some modular group $\Gamma(N) \subset SL(2, \Z)$. In physics, the quantum field theory given by this class of representations of $\hat{\mathfrak{g}}$ is called an admissible level WZW model.

\subsection{The main result of the present paper}

Let $\ell$ be an admissible number. We study a category of (ordinary) modules for the admissible affine vertex operator algebra $L(\ell, 0)$. Recall that a (ordinary) module by definition satisfies two grading-restriction conditions---finite dimensionality of the weight spaces and lower boundedness of the conformal weights. We denote this category by $\mathcal{O}_{\ell, \ord}$. It is a subcategory of the category $\mathcal{O}$ (denoted by $\mathcal{O}_{\ell}$ in this paper)  studied by Arakawa in \cite{A} and is proved in \cite{A} to be semisimple and to have finitely many irreducible objects up to equivalence. In this paper, we prove that  the category $\mathcal{O}_{\ell, \ord}$ has a natural braided tensor category structure (with a twist)
by using the logarithmic tensor category theory developed in \cite{HLZ0}--\cite{HLZ8} and \cite{H4}. The work carried out in the present paper is the verification that the assumptions in these papers are indeed satisfied by 
the category $\mathcal{O}_{\ell, \ord}$.  

It is easy to see that an $L(\ell, 0)$-module in $\mathcal{O}_{\ell, \ord}$ is $C_{1}$-cofinite. By a result of the second author in 
\cite{H3}, the fusion rules among modules in $\mathcal{O}_{\ell, \ord}$ are all finite. This finiteness together with the results of Arakawa \cite{A} shows that $\mathcal{O}_{\ell, \ord}$  is finitely reductive, in the sense that it is semisimple, has finitely many irreducible objects up to equivalence and the fusion rules among objects in $\mathcal{O}_{\ell, \ord}$ are finite.  As a consequence, the $P(z)$-tensor product, defined in \cite{HLZ3} (cf. \cite{HL2} and \cite{HL3}), exists in $\mathcal{O}_{\ell, \ord}$.

Note that
though the category $\mathcal{O}_{\ell, \ord}$ is semisimple and has only finitely many irreducible objects, there is a crucial difference between $L(\ell, 0)$ for admissible $\ell$ and $L(\ell, 0)$ for positive integral $\ell$. In the case that $\ell\in \Z_{+}$, every weak module is a direct sum of (ordinary) irreducible modules (in particular, every generalized module is a direct sum of (ordinary) irreducible modules). In the case that $\ell$ is admissible, weak modules or even generalized modules in general might not be completely reducible. This difference
makes our results in this paper much harder to prove than those in the case $\ell\in \Z_{+}$.

Because of the crucial difference discussed above, we cannot use a result of  Frenkel and Zhu \cite{FZ} and Li \cite{Li2}  to calculate or even just 
estimate the fusion rules since this result holds only when every $\N$-gradable weak module is a direct sum of ordinary modules. Instead, we use a result of the second and third authors \cite{HY} to obtain such an estimate. 

 In \cite{H1, HLZ6, HLZ7}, one condition needed in the  construction of  the associativity isomorphism is
that a finitely generated lower-bounded generalized module  is an object in the category. 
But  because of the crucial difference discussed above, in general a finitely generated lower-bounded generalized module (or even a finitely generated $L_0$-semisimple lower-bounded generalized module) might not be an ordinary module. In \cite{H4}, this condition is replaced by a weaker condition  (see Section 4 for the details). Verifying this condition is one of the main work that we carry out in this paper. We need to use some results of \cite{A} and \cite{KL2}. 

Another condition needed  in the  construction of  the associativity isomorphism is the convergence and extension property without logarithm in \cite{H1}. 
Again because of the crucial difference discussed above, it is difficult to prove the convergence and extension property without logarithm directly (see \cite{H3} for a general proof in the case that  every generalized module is a direct sum of (ordinary) irreducible modules). Because of this difficulty, instead of using the semisimple tensor category theory developed in \cite{HL1}--\cite{HL6} and \cite{H1}, we use the logarithmic tensor category developed in \cite{HLZ0}--\cite{HLZ8} and \cite{H4}. In this theory, we need only prove the convergence and extension property with logarithm in \cite{HLZ7}. For $V$-modules satisfying the $C_1$-cofiniteness condition, the second author proved  in \cite{H3} (cf. \cite{HLZ7, Y2}) that matrix elements of products and iterates of intertwining operators satisfy a system of differential equations such that the prescribed singular points are regular. The convergence and extension property with logarithm  follows easily from the properties of these differential equations, as was shown in \cite{H3}. 
Since objects in $\mathcal{O}_{\ell, \ord}$  are all $C_{1}$-cofinite,  the convergence and extension property with logarithm holds. We can also use the Knizhnik-Zamolodchikov differential equations (\cite{KZ, HL7}) to prove this property.

Now since the assumptions in \cite{HLZ0}--\cite{HLZ8} and \cite{H4} all hold, we conclude that $\mathcal{O}_{\ell, \ord}$ has a natural  braided tensor category structure. Moreover it is in fact a semisimple category with finitely many irreducible objects. In particular, every logarithmic intertwining operator among objects in this category is an ordinary intertwining operator without logarithm. Then from the associativity of intertwining operators, we see that in fact the convergence and extension property without logarithm holds.

Having constructed the braided tensor category structure on $\mathcal{O}_{\ell, \ord}$, the next goal is to determine if every object is rigid. 
We answer this question in the special case that the affine Lie algebra is $\widehat{\mathfrak{sl}}_2$.  The admissible levels are $\ell=-2+\frac{a}{b}$ with $a, b$ co-prime positive integers and $a>1$. We prove that every object of $\mathcal{O}_{\ell, \ord}$ is rigid and moreover $\mathcal{O}_{\ell, \ord}$ is a modular tensor category if and only if $b$ is odd.

In this paper, as is mentioned above, we construct a braided tensor category structure on the category of ordinary modules of an affine vertex algebra at an admissible level. 
However, it has already been realized that an affine vertex operator algebra at an admissible level possesses quite complicated modules, i.e. the number of inequivalent irreducible modules is not countable, 
conformal weights of an indecomposable module might not be bounded from below, $L_0$ might not act semisimply and generalized $L_0$-eigenspaces might not be finite dimensional. These features have been mainly explored in the case of affine vertex operator algebras of $\mathfrak{sl}_2$ at an admissible level \cite{G, R1, R2, R3, CR1, CR2, AM, RW}. 
It is a very important problem to develop more techniques in order to study these module categories.
We remark that in the instance of $\mathfrak{sl}_2$ at an admissible negative level, indecomposable modules of its Heisenberg coset $\text{Com}\left(H, L(\ell, 0)\right)$ (here $H$ denotes the Heisenberg vertex operator subalgebra of $L(\ell, 0)$) seem to behave better in the sense that they have lower-bounded conformal weights and finite-dimensional generalized $L_0$-eigenspaces (see \cite{CRW, Ad2, ACR} for examples). Moreover at least in the case of $\ell=-1/2$ and $\ell=-4/3$ all simple coset modules are $C_1$-cofinite \cite[Theorem 13]{CMR} and in general $C_1$-cofiniteness has been proven for a certain subcategory of coset modules \cite[Theorem 6.1]{CKLR}. 
A possible route would thus be to establish braided tensor category structure on coset modules and then use the theory of vertex operator algebra extensions \cite{HKL, CKM} to construct a corresponding category of affine vertex operator algebra modules as the category of local modules of an infinite order simple current extension \cite{CKL}.

\subsection{Conjectures and the case of $\mathfrak{g}=\mathfrak{sl}_2$}

In the case of $\ell+h^{\vee}\not\in \Q_{\ge 0}$ or $\ell\in \Z_{+}$, the most difficult part of the construction is the proof of the rigidity. 
In the case of $\ell+h^{\vee}\not\in \Q_{\ge 0}$, results in the theory of quantum groups are needed (see \cite{KL5}). 
In the case of $\ell\in \Z_{+}$, the Verlinde formula is needed (see \cite{H6} and \cite{F3}). 
We now formulate the rigidity of  $\mathcal O_{\ell, \ord}$ as the first of three conjectures for admissible $\ell$ and
 we prove the first two conjectures in the case of $\mathfrak{g}=\mathfrak{sl}_2$ in Section 7. 

Let $\mathfrak g$ be a simple Lie algebra, $h^\vee$ the dual Coxeter number of $\mathfrak g$ and $d$ its ``lacety'' (which is $1$ for types $A,D,E$, $2$ for $B,C,F$, and $3$ for $G$).
 Let $\ell=-h^\vee+\frac{a}{b}$ with coprime positive integers $a, b$ such that $\ell$ is an admissible level for $\mathfrak g$. 
 Let $q=e^{\frac{2\pi i }{2d(\ell+h^\vee)}}$ and $\mathcal C_\ell(\mathfrak g)$ be the semi-simplification of the category of tilting modules of $U_q(\mathfrak g)$ as described in \cite{Sawin}. Let $\mathcal O_{\ell, \ord}$ be the braided tensor category constructed in the present paper. 
 
The first conjecture is on the rigidity:
\begin{conj}
The category $\mathcal O_{\ell, \ord}$ is ribbon.
\end{conj}
In section 7, we prove this conjecture for $\mathfrak g=\mathfrak{sl}_2$.  

The second conjecture is on the modulartity: 
\begin{conj}
The category $\mathcal O_{\ell, \ord}$ is modular except for the following list ($n\in\mathbb Z_{>0}$)
\begin{enumerate}
\item $\mathfrak g \in \{\mathfrak{sl}_{2n}, \mathfrak{so}_{2n}, \mathfrak{e}_7, \mathfrak{sp}_n\}$ and $b$ even. 
\item $\mathfrak g = \mathfrak{so}_{4n+1}$ and $b=0 \mod 4$.
\item $\mathfrak g = \mathfrak{so}_{4n+3}$ and $b=2 \mod 4$. 
\end{enumerate}
\end{conj}
In section 7, we prove this conjecture for $\mathfrak g=\mathfrak{sl}_2$.  

A stronger conjecture whose correctness would imply the previous two (see \cite[Theorem 6]{Sawin})  is the following:
\begin{conj}
The category $\mathcal O_{\ell, \ord}$ and the semi-simplification $\mathcal C_\ell(\mathfrak g)$ of the category of tilting modules for $U_q(\mathfrak g)$ are equivalent as braided tensor categories. 
\end{conj}

Recently, Arakawa and Kawasetsu \cite{AK} introduced the notion of quasi-lisse vertex algebras. These are strongly finitely generated vertex operator algebras whose associated variety has only finitely many symplectic leaves. Characters of ordinary modules of quasi-lisse vertex algebras converge to holomorphic functions on the upper half of the complex plane and they satisfy a modular differential equation. An affine vertex operator algebra at admissible level is the prime example of a quasi-lisse vertex algebra. Quasi-lisse vertex algebras, their tensor categories and modularity of their characters play an important role in four-dimensional supersymmetric gauge theories \cite{BR, BLL+}. It is thus another interesting problem to find out conditions that ensure that the category of ordinary modules of a quasi-lisse vertex algebra is a (rigid) braided tensor category.

\subsection{Organization}

The paper is organized as follows: In Section 2, we briefly review affine Lie algebras and their modules, vertex operator algebras and their modules associated to affine Lie algebras. In Section 3, we discuss the category $\mathcal{O}_{\ell, \ord}$ for admissible $\ell$. In Section 4, we give an estimate of the fusion rules for admissible affine Lie algebras in the category $\mathcal{O}_{\ell, \ord}$. In Section 5, we recall the definition and construction of the $P(z)$-tensor product and show that the category $\mathcal{O}_{\ell, \ord}$ is closed under the $P(z)$-tensor product operation. In Section 6, we construct the associativity isomorphism and obtain the main result. In Section 7, we prove that when $\mathfrak{g}=\mathfrak{sl}_2$, the category $\mathcal{O}_{\ell, \ord}$ is in fact rigid and even modular in some cases. 

\paragraph{Acknowledgments}

T. C. is supported by the Natural Sciences and Engineering Research Council of Canada (RES0020460). J. Y. is supported in part by an AMS-Simons travel grant. J. Y. also wants to thank Zongzhu Lin for useful conversations.

\setcounter{equation}{0}
\section{Vertex operator algebras associated with affine Lie algebras}
\subsection{Affine Lie algebras and their modules}

Let $\mathfrak{g}$ be a complex simple Lie algebra of rank $r$. Fix a triangular decomposition
\[
\mathfrak{g} = \mathfrak{n}_{-} \oplus \mathfrak{h} \oplus \mathfrak{n}_{+},
\]
with a Cartan subalgebra $\mathfrak{h}$ of $\mathfrak{g}$. We will often identify $\mathfrak{h}$ with $\mathfrak{h}^{*}$ using the invariant symmetric bilinear form $(\cdot, \cdot)$.

For $\lambda \in \mathfrak{h}^{*}$, we use $L(\lambda)$ to denote the irreducible highest weight $\mathfrak{g}$-module with highest weight $\lambda$.

The affine Lie algebra $\hat{\mathfrak{g}}$ associated with $\mathfrak{g}$ and $(\cdot, \cdot)$ is the vector space $\mathfrak{g} \otimes \C[t, t^{-1}] \oplus \C{\bf k}$ equipped with the bracket operation defined by
\[
[a \otimes t^m, b \otimes t^n] = [a, b]\otimes t^{m+n} + (a, b)m\delta_{m+n,0}{\bf k},
\]
\[
[a \otimes t^m, {\bf k}] = 0,
\]
for $a, b \in \mathfrak{g}$ and $m, n \in \Z$. It is $\Z$-graded in a natural way. Consider the subalgebras
\[
\hat{\mathfrak{g}}_{\pm} = \mathfrak{g}\otimes t^{\pm 1}\C[t^{\pm 1}]
\]
and the vector space decomposition
\[
\hat{\mathfrak{g}} = \hat{\mathfrak{g}}_{-} \oplus \mathfrak{g} \oplus \C{\bf k} \oplus \hat{\mathfrak{g}}_{+}.
\]

Let $M$ be a $\mathfrak{g}$-module, viewed as homogeneously graded of fixed degree, and let $\ell \in \C$. Let $\hat{\mathfrak{g}}_{+}$ act on $M$ trivially and ${\bf k}$ act as the scalar multiplication by $\ell$. Then $M$ becomes a $\mathfrak{g} \oplus \C{\bf k} \oplus \hat{\mathfrak{g}}_{+}$-module, and we have the $\C$-graded induced $\hat{\mathfrak{g}}$-module
\[
\widehat{M}_{\ell} = U(\hat{\mathfrak{g}})\otimes_{U(\mathfrak{g}\oplus \C{\bf k}\oplus \hat{\mathfrak{g}}_{+})}M,
\]
which contains a canonical copy of $M$ as its subspace.

For $\lambda \in \mathfrak{h}^{*}$, we use the notation $M(\ell, \lambda)$ to denote the graded $\hat{\mathfrak{g}}$-module $\widehat{L(\lambda)}_{\ell}$. Let $J(\ell, \lambda)$ be the maximal proper submodule of $M(\ell, \lambda)$ and $L(\ell, \lambda) = M(\ell, \lambda)/J(\ell, \lambda)$. Then $L(\ell, \lambda)$ is the unique irreducible graded $\hat{\mathfrak{g}}$-module such that ${\bf k}$ acts as $\ell$ and the space of all elements annihilated by $\hat{\mathfrak{g}}_{+}$ is isomorphic to the $\mathfrak{g}$-module $L(\lambda)$.

\subsection{Vertex operator algebras associated with affine Lie algebras}
In the special case $\lambda = 0$, $L(0)$ is the one-dimensional  trivial $\mathfrak{g}$-module. Consequently, as a vector space and $\hat{\mathfrak{g}}_{-}$-module,
\[
M(\ell, 0) = U(\hat{\mathfrak{g}})\otimes_{U(\mathfrak{g}\oplus \C{\bf k}\oplus \hat{\mathfrak{g}}_{+})}\C \simeq U(\hat{\mathfrak{g}}_{-})\otimes \C = U(\hat{\mathfrak{g}}_{-}).
\]
Note that $M(\ell, 0)$ is spanned by the elements of the form $a_1(-n_1)\cdots a_m(-n_m)1$, where $a_1, \dots, a_m \in \mathfrak{g}$ and $n_1, \dots, n_m \in \Z_{+}$, with $a(-n)$ denoting the representation image of $a \otimes t^{-n}$ for $a \in \mathfrak{g}$ and $n \in \Z_{+}$. We define the vertex operator map recursively on $m$ as follows: For $1 \in M(\ell, 0)$, we define $Y(1, x)$ to be the identity operator on $M(\ell, 0)$ and for $a \in \mathfrak{g}$, we define
\[
Y(a(-1)1, x) = \sum_{n \in \Z}a(n)x^{-n-1}.
\]
Assume that the map has been defined for $a_1(-n_1)\cdots a_m(-n_m)1$. Vertex operator for elements of the form $a_0(-n_0)a_1(-n_1)\cdots a_m(-n_m)1$ are defined using the residue in $x_1$ of the Jacobi identity for vertex operator algebras as follows:
\begin{eqnarray}\label{eq11}
&&Y(a_0(-n_0)a_1(-n_1)\cdots a_m(-n_m)1, x)\nn
&=& \res_{x_1}(x_1-x)^{-n_0}Y(a_0(-1)1, x_1)Y(a_1(-n_1)\cdots a_m(-n_m)1, x)\nn
&& - \res_{x_1}(-x+x_1)^{-n_0}Y(a_1(-n_1)\cdots a_m(-n_m)1, x)Y(a_0(-1)1, x_1).\nn
\end{eqnarray}
Here and below, binomial expressions are understood to be expanded in nonnegative powers of the second variable. The vacuum vector for $M(\ell, 0)$ is ${\bf 1} = 1$. In the case that $\ell \neq -h^{\vee}$, $M(\ell,0)$ also has a Virasoro element
\begin{equation}\label{eq12}
\omega = \frac{1}{2(k+h^{\vee})}\sum_{i=1}^{\mbox{dim}\;\mathfrak{g}}g^i(-1)^2{\bf 1},
\end{equation}
where $\{g^i\}_{i=1, \dots, \mbox{dim}\; \mathfrak{g}}$ is an arbitrary orthonormal basis of $\mathfrak{g}$ with respect to the form $(\cdot, \cdot)$. With respect to the grading, ${\bf 1}$ has (conformal) weight $0$; the weight of a homogeneous element is the negative of its degree. The generalized Verma module $M(\ell, 0)$ has a vertex operator structure:
\begin{thm}[\cite{FZ}]\label{FZ1} If $\ell \neq -h^{\vee}$, the quadruple $(M(\ell,0), Y, {\bf 1}, \omega)$ defined above is a vertex operator algebra.
\end{thm}

Since $J(\ell,0)$ is a $\hat{\mathfrak{g}}$-submodule of $M(\ell,0)$, the vertex operator map for $M(\ell,0)$ induces a vertex operator map for $L(\ell,0)$ which we still denote by $Y$. We continue to denote the $J(\ell,0)$-cosets of ${\bf 1}$ and $\omega$ in $M(\ell,0)$ by ${\bf 1}$ and $\omega$. Thus we have the following consequence of Theorem \ref{FZ1}:
\begin{cor}[\cite{FZ}]If $\ell \neq -h^{\vee}$, the quadruple $(L(\ell,0), Y, {\bf 1}, \omega)$ defined above is a vertex operator algebra.
\end{cor}

For $\ell \neq -h^{\vee}$ and $\lambda \in \mathfrak{h}^*$, the generalized Verma module $M(\ell, \lambda)$ is a module for $M(\ell,0)$ with the vertex operator map
\begin{equation}\label{eq13}
Y(\cdot, x): M(\ell,0) \longrightarrow (\edo M(\ell, \lambda))[[x, x^{-1}]]
\end{equation}
defined recursively just as in (\ref{eq11}). For $\ell \in \Z_{+}$, the following theorem was proved in \cite{FZ}:
\begin{thm}[\cite{FZ}]
Let $\ell \in \Z_{+}$. Then every $\N$-gradable weak $L(\ell,0)$-module is a direct sum of irreducible $L(\ell,0)$-module  and
\[
\{L(\ell, \lambda)|\lambda \in \mathfrak{h}^{*}\; \mbox{is dominant integral such that}\; (\lambda, \theta) \leq \ell\}
\]
is the set of all irreducible $L(\ell, 0)$-modules up to equivalence, with the vertex operator maps induced from (\ref{eq13}).
\end{thm}

\section{The category $\mathcal{O}_{\ell, \ord}$ for an admissible level $\ell$}
In the remaining of this paper, unless further stated, we will let $\ell$ be an admissible number and study the category of $L(\ell, 0)$-modules.

Let $\mathcal{O}_{\ell}$ be the full subcategory of $\hat{\mathfrak{g}}$-modules of level $\ell$ satisfying the following conditions:
\begin{enumerate}
\item[(1)]$L_0$ acts locally finitely;
\item[(2)]$\hat{\mathfrak{g}}_{+} \oplus \mathfrak{n}_{+}$ acts locally nilpotently;
\item[(3)]$\hat{\mathfrak{h}}: = \mathfrak{h} \oplus \C{\bf k}$ acts semisimply.
\end{enumerate}

We recall that 
a weak module for a vertex operator algebra is a vector space equipped with a vertex operator map satisfying all the axioms for a module except those axioms involving the grading and the action of $L_0$. We also recall that a generalized module for a vertex operator algebra is weak module equipped with a $\C$-grading on the underlying vector space such that 
the homogeneous subspaces are generalized eigenvectors of the operator $L_{0}$. 

The following theorem is the main theorem in \cite{A}:
\begin{thm}[\cite{A}]\label{admmain}
Let $\ell$ be an admissible number. Then up to equivalence, there are only finitely many irreducible $L(\ell,0)$-modules in the category $\mathcal{O}_{\ell}$.
Moreover, every generalized $L(\ell,0)$-module in the category $\mathcal{O}_{\ell}$ is a direct sum of irreducible $L(\ell,0)$-modules in $\mathcal{O}_{\ell}$.
\end{thm}

\begin{rema}
{\rm Note that every object in $\mathcal{O}_{\ell}$ has a grading given by the generalized eigenspaces of $L_{0}$. In particular, 
every weak $L(\ell,0)$-module in the category $\mathcal{O}_{\ell}$ is in fact a generalized $L(\ell,0)$-module. Hence by Theorem \ref{admmain},
every weak $L(\ell,0)$-module in the category $\mathcal{O}_{\ell}$  is a direct sum of irreducible $L(\ell,0)$-modules in $\mathcal{O}_{\ell}$.}
\end{rema}

Let $\mathcal{O}_{\ell, \ord}$ be the category of $L(\ell,0)$-modules. In particular, an element in $\mathcal{O}_{\ell, \ord}$ possesses a decomposition of $L_{0}$-weight spaces satisfying the two grading restrictions: the $L_{0}$-weights are lower bounded and the $L_{0}$-weight spaces are finite dimensional. 

\begin{cor}\label{admlemmamain}
Let  $\ell$ be an admissible number. Then
\begin{itemize}
\item[(1)] The category $\mathcal{O}_{\ell, \ord}$ is a subcategory of $\mathcal{O}_{\ell}$. In particular, 
there are only finitely many irreducible $L(\ell,0)$-modules, up to equivalence.
\item[(2)]Every $L(\ell,0)$-module is completely reducible, that is, the category $\mathcal{O}_{\ell, \ord}$ is semisimple.
\end{itemize}
\end{cor}
\pf 
Let $W$ be a simple object in $\mathcal{O}_{\ell, \ord}$. Then the lowest weight subspace of $W$, which is finite dimensional, must be of the form $L(\lambda)$ for some dominant integral weight $\lambda \in \mathfrak{h}^*$. Thus $W \simeq L(\ell, \lambda)$. From \cite{A}, $\lambda$ must be an admissible weight and from \cite{GK} (see also \cite{A}), an object in $\mathcal{O}_{\ell, \ord}$ is a direct sum of simple objects in $\mathcal{O}_{\ell, \ord}$. \epfv

We denote by $P$ the set of weights $\lambda \in \mathfrak{h}^*$ such that $L(\ell, \lambda)$ are the simple objects in $\mathcal{O}_{\ell, \ord}$. In particular, $P$ is a finite set.

We recall the following notations in \cite{KL2}: For any integer $N \geq 1$, let $Q_N$ to be the subspace of $U(\hat{\mathfrak{g}})$ spanned by the products $(a_1 \otimes t)\dots (a_N \otimes t)$ with $a_1, \dots, a_N \in \mathfrak{g}$. (Here and below we use the notation $U(\mathcal{L})$ for the universal enveloping algebra of a Lie algebra $\mathcal{L}$.) For a $\hat{\mathfrak{g}}$-module $W$ and an integer $N \geq 1$, let
\[
W(N) = \{w \in W\;| \; Q_Nw = 0\}.
\]
It is easy to see that
\[
W(1) \subset W(2) \subset \cdots.
\]
Let
\[
W(\infty) = \bigcup_{N \geq 1}W(N).
\]
Recall from \cite{KL2} that $W$ is {\it smooth} if $W = W(\infty)$.

We also recall that 
a {\em weak module} for a vertex operator algebra is a vector space equipped with a vertex operator map satisfying all the axioms for a module except those axioms involving the grading and the action of $L_0$.

We now give a characterization of the category $\mathcal{O}_{\ell, \ord}$. 
\begin{thm}\label{charcatO}
Let $\ell$ be an admissible level and let $W$ be a weak $L(\ell,0)$-module. Then $W$ is in $\mathcal{O}_{\ell, \ord}$ if and only if $W$ viewed as a 
$\hat{\mathfrak{g}}$-module is smooth and dim $W(1) < \infty$.
\end{thm}
\pf The ``only if" part is obvious from Theorem \ref{admmain} and Corollary \ref{admlemmamain}. We will prove the ``if" part. 

Let $W$ be a weak $L(\ell,0)$-module. Assume that $W$ viewed as a 
$\hat{\mathfrak{g}}$-module is smooth and dim $W(1) < \infty$. For any integer $N \geq 1$, from the definition of the subspace $W(N)$,  it is a 
$\mathfrak{g}$-module. From Lemma 1.10 in \cite{KL2}, $W(N)$ is also finite dimensional. (Here and below we note that in \cite{KL2}, 
the restriction on the level is assumed starting from Definition 2.15.) 
In particular, $\mathfrak{h}$ acts semisimply and $\mathfrak{n}_{+}$ acts locally nilpotently on $W(N)$.
 
Let $W\{N\}$ be the weak $L(\ell,0)$-submodule of $W$ generated by $W(N)$. From Lemma 2.5 of \cite{KL2}, $W\{N\}$ is a quotient of a generalized Weyl module induced from $W(N)$. Thus $W\{N\}$ is a generalized $L(\ell,0)$-module satisfying the two grading restriction conditions. In particular, $L_{0}$ acts locally finitely and $\mathfrak{g}_{+}$ acts locally nilpotently on $W\{N\}$.
Also, because $\mathfrak{h}$ acts semisimply and $\mathfrak{n}_{+}$ acts locally nilpotently on $W(N)$, $\hat{\mathfrak{h}}$ acts semisimply and $\mathfrak{n}_{+}$ acts locally nilpotently on $W\{N\}$. Thus $W\{N\}$ is in $\mathcal{O}_{\ell}$. Since $W\{N\}$ is a generalized $L(\ell,0)$-module satisfying the grading-restriction conditions,  by Threorem \ref{admmain}, $W\{N\}$ is in $\mathcal{O}_{\ell, \ord}$.

Since $W\{N\}$ is in $\mathcal{O}_{\ell, \ord}$, the length of $W\{N\}$ is the number of irreducible components of $W\{N\}(1)$. Since $W\{N\}$ is a submodule of $W$, we have $W\{N\}(1) \subset W(1)$.  In particular, the length of $W\{N\}$ is bounded above by the number of irreducible components in the finite-dimensional $\mathfrak{g}$-module  $W(1)$. Since this this upper bound of the lengths of $W\{N\}$ for $N\ge 1$ is independent of $N$, the ascending sequence of $L(\ell,0)$-modules $W\{1\} \subset W\{2\} \subset \cdots$ must be stationary. Since $W$ is smooth, the union of $W\{N\}$ for all $N \geq 1$ is equal to $W$. It follows that $W = W\{N\}$ for some $N \geq 1$. Since $W\{N\}$ is in $\mathcal{O}_{\ell, \ord}$, the same must hold for $W$. \epfv

\section{Intertwining operators and fusion rules}

\begin{defn}\label{log:def}{\rm
Let $(W_1,Y_1)$, $(W_2,Y_2)$ and $(W_3,Y_3)$ be generalized modules
for a vertex operator algebra $V$. A {\em logarithmic intertwining
operator of type ${W_3\choose W_1\,W_2}$} is a linear map
\begin{align}\label{log:map0}
\mathcal{Y}: &W_1\otimes W_2\to W_3\{x\}[\log x]\nn
&w_{(1)}\otimes w_{(2)}\mapsto{\mathcal{Y}}(w_{(1)},x)w_{(2)}=\sum_{k=0}^{K}\sum_{n\in
{\mathbb C}}{w_{(1)}}_{n, k}^{
\mathcal{Y}}w_{(2)}x^{-n-1}(\log x)^{k}
\end{align}
for all $w_{(1)}\in W_1$ and $w_{(2)}\in W_2$, such that the
following conditions are satisfied:
(i) The {\em lower truncation
condition}: for any $w_{(1)}\in W_1$, $w_{(2)}\in W_2$, $n\in
{\mathbb C}$ and $k=0, \dots, K$,
\begin{equation}
{w_{(1)}}_{n+m, k}^{\mathcal{Y}}w_{(2)}=0\;\;\mbox{ for }\;m\in {\mathbb
N} \;\mbox{sufficiently large}.
\end{equation}
(ii) The {\em Jacobi identity}:
\begin{eqnarray}\label{Jacob}
\lefteqn{\dps x^{-1}_0\delta \bigg( {x_1-x_2\over x_0}\bigg)
Y_3(v,x_1){\mathcal{Y}}(w_{(1)},x_2)w_{(2)}}\nn
&&\hspace{2em}- x^{-1}_0\delta \bigg( {x_2-x_1\over -x_0}\bigg)
{\mathcal{Y}}(w_{(1)},x_2)Y_2(v,x_1)w_{(2)}\nn
&&{\dps = x^{-1}_2\delta \bigg( {x_1-x_0\over x_2}\bigg){
\mathcal{Y}}(Y_1(v,x_0)w_{(1)},x_2) w_{(2)}}
\end{eqnarray}
for $v\in V$, $w_{(1)}\in W_1$ and $w_{(2)}\in W_2$. (iii) The {\em $L_{-1}$-derivative property:} for any
$w_{(1)}\in W_1$,
\begin{equation}
{\mathcal{Y}}(L_{-1}w_{(1)},x)=\frac d{dx}{\mathcal{Y}}(w_{(1)},x).
\end{equation}
A logarithmic intertwining operator is called an {\it intertwining operator} if no $\log x$ appears in ${\mathcal{Y}}(w_{(1)},x)w_{(2)}$
for $w_{(1)}\in W_1$ and $w_{(2)}\in W_2$. The dimension of the space  $\mathcal{V}_{W_1 W_2}^{W_3}$ of all logarithmic intertwining 
operators is called the {\it fusion rule} of the same type and is denoted $\mathcal{N}_{W_1 W_2}^{W_3}$.}
\end{defn}

\begin{rema}
{\rm When $(W_1,Y_1)$, $(W_2,Y_2)$ and $(W_3,Y_3)$ are $V$-modules, a logarithmic 
intertwining operator of type ${W_3\choose W_1\,W_2}$ must
be an intertwining operator of the same type.}
\end{rema}

We now discuss fusion rules. It is clear that every irreducible $L(\ell, 0)$-module in $\mathcal{O}_{\ell, \ord}$
is $C_{1}$-cofinite and the contragredient module of an irreducible $L(\ell, 0)$-module in $\mathcal{O}_{\ell, \ord}$
is still in $\mathcal{O}_{\ell, \ord}$. Thus by Theorem 3.1 in \cite{H3}, the fusion rules among these irreducible
$L(\ell, 0)$-modules are all finite. This finiteness and Theorem \ref{admmain} can be summarized together to 
give the following:
\begin{thm}\label{fin-red}
If $\ell$ is admissible, the vertex operator algebra $L(\ell,0)$ is finitely reductive in the sense that
\begin{itemize}
\item[(1)]The category $\mathcal{O}_{\ell, \ord}$ is a subcategory of $\mathcal{O}_{\ell}$. In particular, 
there are only finitely many irreducible $L(\ell,0)$-modules, up to equivalence.
\item[(2)]Every $L(\ell,0)$-module is completely reducible, that is, the category $\mathcal{O}_{\ell, \ord}$ is semisimple.
\item[(3)]The fusion rules for $L(\ell, 0)$ are finite.
\end{itemize}
\end{thm}

This finiteness of the fusion rules is enough for the main results in this paper. But here we give a 
stronger result than this finiteness by giving an estimate of the fusion rules in terms of the corresponding
irreducible $\mathfrak{g}$-modules. 

In \cite{FZ}, I. Frenkel and Zhu identified the spaces of intertwining operators among irreducible 
ordinary $V$-modules with homomorphisms between certain
modules for Zhu's algebra $A(V)$. In fact, this result holds only in the case that every 
$\N$-gradable weak $V$-module is a  direct sums of irreducible $V$-modules. A proof of this result
and a discussion on the importance of this condition are given by Li in \cite{Li2}. 
For the vertex operator algebra $L(\ell, 0)$ that we are discussing in this paper, not 
every $\N$-gradable weak $V$-module is a  direct sums of irreducible $V$-modules. In particular, we cannot use this 
result to calculate or even estimate the fusion rules. Instead, we shall use a result in \cite{HY} to give an estimate
of the fusion rules among irreduible $L(\ell, 0)$-modules in $\mathcal{O}_{\ell, \ord}$.

We recall the definition of Zhu's algebra for a vertex operator algebra $V$ from \cite{Z}. 
The vertex operator algebra $V$ has a product $*$ given by
\[
u*v = \res_xx^{-1}Y((1+x)^{L_0}u, x)v
\]
for $u, v \in V$. We also define the subspace $O(V) \subset V$ as the linear span of elements
\[
\res_xx^{-2}Y((1+x)^{L_0}u,x)v
\]
for $u, v \in V$. The $*$ is well-defined product on the quotient $A(V) = V/O(V)$, and in fact $(A(V), *)$ is an associative algebra with unit ${\bf 1} + O(V)$, called the {\em Zhu's algebra} of $V$.

If $W$ is a $V$-module, then there is a left action of $V$ on $W$ defined by
\[
u*w = \res_xx^{-1}Y((1+x)^{L_0}u, x)w
\]
for $u \in V$, $w \in W$, and a right action defined by
\[
w*v = \res_xx^{-1}Y((1+x)^{L_0-1}v, x)w
\]
for $v \in V$, $w \in W$. If we define the subspace $O(W) \subset W$ as the linear span of elements of the form
\[
\res_xx^{-2}Y((1+x)^{L_0}u,x)w
\]
for $u$ and $w$, then the left and right actions define an $A(V)$-bimodule structure on the quotient $A(W) = W/O(W)$.

For a $V$-module $W$ we define the {\em top level} $T(W)$ to be the subspace
\[
T(W) = \{w \in W|v_nw = 0, v \in V \; \mbox{homogeneous},\; \wt v - n -1 < 0\}.
\]
Then $T(W)$ is an $A(V)$-module with action defined by
\[
(v + O(V))\cdot w  = o(v)w
\]
for $v \in V$ and $w \in T(W)$; here $o(v) = v_{\wt v - 1}$ if $v$ is homogeneous and we extend linearly to define $o(v)$ for non-homogeneous $v$.

The following result from \cite{HY} gives the upper bound for the fusion rules among $V$-modules generated by their top levels:
\begin{prop}[\cite{HY}]\label{fr}
Let $W_1$, $W_2$ and $W_3$ be $V$-modules such that $W_2$ and $W_3'$ are generated by $T(W_2)$ and $T(W_3')$, respectively. Then the map
\begin{eqnarray*}
\rho: \mathcal{V}_{W_1 W_2}^{W_3} &\longrightarrow & \hom_{A(V)}(A(W_1)\otimes_{A(V)}T(W_2), T(W_3))\nn
\mathcal{Y} &\longmapsto & \rho(\mathcal{Y})
\end{eqnarray*}
defined by
\[
\rho(\mathcal{Y})(w_1 \otimes w_2) = \res_x x^{-1}\mathcal{Y}(w_1, x)w_2
\]
for $w_1 \in A(W_1)$ and $w_2 \in T(W_2)$, is well-defined and injective.
\end{prop}

We shall apply Proposition \ref{fr} to the vertex operator algebra $L(\ell,0)$ of admissible level $\ell$ and show that the fusion rule for a triple of modules in $\mathcal{O}_{\ell, \ord}$ is finite. Since every module in the category $\mathcal{O}_{\ell,\ord}$ is completely reducible, we  need only consider the case that 
the $L(\ell,0)$-modules are of the form $W= L(\ell, \lambda)$ such that $\lambda \in P$.

\begin{prop}
Let $W_i = L(\ell, \lambda_i)$ such that $\lambda_i \in P$ for $i = 1, 2, 3$. Then we have the following estimate of the fusion rule:
\begin{equation}\label{fusion-rules}
\mathcal{N}_{W_1 W_2}^{W_3} \leq \mbox{dim}\; \hom_{\mathfrak{g}}(L(\lambda_1)\otimes L(\lambda_2), L(\lambda_3)) < \infty.
\end{equation}
\end{prop}
\pf
From \cite{FZ}, we have
\[
A(L(\ell, 0)) \cong U(\mathfrak{g})/I_{\ell}
\]
where $I_{\ell}$ is  a two sided-ideal of the universal enveloping algebra $U(\mathfrak{g})$ of $\mathfrak{g}$. In general, the ideal $I_{\ell}$ is hard to determine explicitly. Similarly, we have
\begin{equation}\label{zhua}
A(L(\ell, \lambda)) \cong (L(\lambda) \otimes U(\mathfrak{g}))/J_{\ell, \lambda},
\end{equation}
where the left action of $A(L(\ell, 0))$ on $A(L(\ell, \lambda))$ is induced from
\begin{equation}\label{zhub}
a(v \otimes x) = a \cdot v \otimes x + v \otimes ax
\end{equation}
and the right action is induced from
\[
(v \otimes x)a = v \otimes xa
\]
for $a \in \mathfrak{g}$, $v \in L(\lambda)$, $x \in U(\mathfrak{g})$, and $J_{\ell, \lambda}$ is a two-sided $U(\mathfrak{g})/I_{\ell}$-submodule of $L(\lambda) \otimes U(\mathfrak{g})$. From Proposition \ref{fr},
\[
\mathcal{V}_{L(\ell, \lambda_1) L(\ell, \lambda_2)}^{L(\ell, \lambda_3)} \subset \hom_{U(\mathfrak{g})/I_{\ell}}((L(\lambda_1) \otimes U(\mathfrak{g}))/J_{\ell, \lambda} \otimes_{U(\mathfrak{g})/I_{\ell}} L(\lambda_2), L(\lambda_3)).
\]
Thus we obtain (\ref{fusion-rules}).
\epfv

\setcounter{equation}{0}
\section{Tensor products of modules in $\mathcal{O}_{\ell, \ord}$}

In this section, we construct ($P(z)$-)tensor products of modules in $\mathcal{O}_{\ell, \ord}$ for an admissible level $\ell$. 
We first recall a result (Theorem \ref{closedness} below) in \cite{KL2} (cf. \cite{Zh}). Although Kazhdan and Lusztig mostly studied module categories of levels that are not positive rational numbers, this result works for objects in $\mathcal{O}_{\ell, \ord}$ for an admissible $\ell$. We shall need the this result later. 

Let $p_0, p_1, \dots, p_m$ be distinct points (or {\em punctures}) on the Riemann sphere $C = \C\cup \{\infty\}$ and let $\varphi_s: C \rightarrow \C$ be the standard local coordinate maps vanishing at $p_s$ for each $s$ (see \cite{H}). We still use $C$ to denote the Riemann sphere equipped with these punctures and local coordinates.

Let $R$ denote the algebra of regular functions on $C\setminus \{p_0, p_1, \dots, p_n\}$. Define
\[
\Gamma_R = \mathfrak{g} \otimes R \oplus \C{\bf k}
\]
with central element ${\bf k}$ and bracket relations
\[
[a \otimes f_1, b \otimes f_2] = [a, b]\otimes f_1f_2 + \{f_1, f_2\}(a,b){\bf k},
\]
for $a, b \in \mathfrak{g}$ and $f_1, f_2 \in R$ and
\[
\{f_1, f_2\} = \res_{p_0}f_2df_1.
\]
Let $W_1, W_2, \dots, W_m$ be $\hat{\mathfrak{g}}$-modules of level $\ell$ (Here we need only $\ell \neq - h^{\vee}$). Then there is a $\Gamma_R$-module structure on $W_1 \otimes \cdots \otimes W_m$, where ${\bf k}$ acts as the scalar $-\ell$ (see \cite{KL2} and \cite{Zh} for details). The dual vector space
\[
(W_1 \otimes \cdots \otimes W_m)^* = \hom(W_1 \otimes \cdots \otimes W_m, \C)
\]
has an induced natural $\Gamma_R$-module structure given by
\[
\langle \xi\cdot\lambda, w\rangle = - \langle \lambda, \xi\cdot w\rangle
\]
for all $\xi \in \Gamma_R$, $\lambda \in (W_1 \otimes \cdots \otimes W_m)^*$, and $w \in W_1 \otimes \cdots \otimes W_m$.

Let $N$ be a positive integer. Let $G_N$ be the subspace of $U(\Gamma_R)$ spanned by all products of the form $(a_1 \otimes f_1)\cdots (a_N \otimes f_N)$ with $a_1, \dots, a_N \in \mathfrak{g}$ and $f_1, \dots, f_N \in R$ satisfying $\iota_{p_0}f_i \in t\C[[t]]$ for $i = 1, \dots, N$, where $\iota_{p_0}$ is the operation that expands an element of $R$ as a power series near $p_0$. Define
\[
Z^N = \{\lambda \in (W_1 \otimes \cdots \otimes W_m)^*|G_N\cdot \lambda = 0\}.
\]
Then we have an ascending sequence $Z^1 \subset Z^2 \subset \cdots$. Let
\[
Z^{\infty} = \bigcup_{N \in \Z_{+}}Z^N.
\]
It is clear that $Z^{\infty}$ is a $\Gamma_R$-submodule of $(W_1 \otimes \cdots \otimes W_m)^*$.

Define a $\mathfrak{g}\otimes \C((t)) \oplus \C{\bf k}$-module structure on $Z^{\infty}$ as follows: ${\bf k}$ acts as scalar $\ell$; and for $\lambda \in Z^{\infty}$, $a \in \mathfrak{g}$ and $g \in \C((t))$, choose $N \in \N$ such that $\lambda \in Z^N$, choose $f \in R$ such that $\iota_{p_0}f -g \in t^N\C[[t]]$, and define
\[
(a \otimes g)\cdot \lambda = (a \otimes f)\cdot \lambda.
\]
Restricted to the Lie subalgebra $\hat{\mathfrak{g}}$ of $\mathfrak{g} \otimes \C((t)) \oplus \C{\bf k}$, we have on $Z^{\infty}$ a structure of $\hat{\mathfrak{g}}$-module of level $\ell$. 

Now we consider the case $m = 2$. We denote the module $Z^{\infty}$ by $W_1 \circ_C W_2$ using the notation in \cite{Zh} (In \cite{KL2}, $Z^{\infty}$ is denoted by $T'(W_1 \otimes W_2)$).

Using Corollary 7.5 in  \cite{KL2} and Theorem \ref{charcatO}, we obatin:
\begin{thm}\label{closedness}
If $W_1, W_2 \in \mathcal{O}_{\ell, \ord}$, then $W_1 \circ_C W_2 \in \mathcal{O}_{\ell, \ord}$.
\end{thm}
\pf From the second paragraph in brackets on page 938 in \cite{KL2}, the proof of  Corollary $7.5$ in \cite{KL2} does not use the assumption there on $\ell$. In particular, we can apply Corollary $7.5$ in \cite{KL2}.  By this corollary, $W_1 \otimes W_2/G_1\cdot(W_1 \otimes W_2)$ is finite dimensional.

By definition, $W_1 \circ_C W_2$ is a smooth $\hat{\mathfrak{g}}$-module. Since
\[
W_1 \circ_C W_2(1) = Z^1 = \hom_{\C}(W_1 \otimes W_2/G_1\cdot(W_1 \otimes W_2), \C)
\]
is finite dimensional, by Theorem \ref{charcatO}, $W_1 \circ_C W_2 \in \mathcal{O}_{\ell, \ord}$. \epfv

Let $z$ be a nonzero complex number. Consider the Riemann sphere $C$ with punctures $p_0 = \infty, p_1 = z$ and $p_2 = 0$ and local coordinate maps vanishing at $p_0, p_1, p_2$ given by $\varphi_0(\epsilon) =1/\epsilon $, $\varphi_1(\epsilon) = \epsilon - z$,  $\varphi_2(\epsilon) = \epsilon$, respectively. A Riemann sphere equipped with these punctures and local coordinates is denote by $P(z)$ as in \cite{HL3}. Then by definition, 
\begin{eqnarray}\label{defnKL}
\lefteqn{W_1 \circ_{P(z)} W_2}\nn
&&=  \{\lambda \in (W_1 \otimes W_2)^*\;|\;\mbox{for some}\; N \in \N,\; (\xi_1\cdots \xi_N)\cdot \lambda = 0,\nn
&& \quad\quad\quad\quad \;\;\;\mbox{for any}\; \xi_1,\dots, \xi_N \in \mathfrak{g}\otimes t\C[t, (z+t)^{-1}]\}.
\end{eqnarray}

As is in \cite{KL2}, we will see later that the contragredient module of $W_1 \circ_{P(z)} W_2$ is the $P(z)$-trensor product module 
of $W_{1}$ and $W_{2}$ that we shall construct later. But note that only a tensor product module
is not enough. We also need addition structures associated with the tensor product module in the construction of the tensor category structure.

Now we recall the construction of $P(z)$-tensor products in \cite{HL3} and \cite{HL5} (see also \cite{HLZ3}--\cite{HLZ4}).

Let $V$ be a vertex operator algebra and let $(W, Y_W)$ be a grading-restricted  $V$-module with
\[
W = \coprod_{n \in \C}W_{[n]}
\]
where for $n \in \C$, $W_{[n]}$ is the generalized weight space with weight $n$. Its {\em contragredient module} is the vector space
\[
W' = \coprod_{n \in \C}(W_{[n]})^*,
\]
equipped with the vertex operator map $Y'$ defined by
\[
\langle Y'(v,x)w', w \rangle = \langle w', Y^{\circ}_W(v,x)w \rangle
\]
for any $v \in V$, $w' \in W'$ and $w \in W$, where
\[
Y^{\circ}_W(v,x) = Y_W(e^{xL(1)}(-x^{-2})^{L_0}v, x^{-1}),
\]
for any $v \in V$, is the {\em opposite vertex operator} (cf. \cite{FHL}).
We also use the standard notation
\[
\overline{W} = \prod_{n \in \C}W_{[n]},
\]
for the formal completion of $W$ with respect to the $\C$-grading.

Fix a nonzero complex number $z$. The concept of $P(z)$-intertwining map, $P(z)$-product and $P(z)$-tensor product are defined as follows 
(see \cite{HL3} and \cite{HLZ3}):

\begin{defn}{\rm
Let $W_1$, $W_2$ and $W_3$ be modules for a vertex operator algebra $V$. A {\em $P(z)$-intertwining map of type ${W_3\choose W_1\,W_2}$} is a linear map
\[
I: W_1 \otimes W_2 \longrightarrow \overline{W_3},
\]
satisfying the following conditions: the {\em lower truncation condition}: for any element $w_{(1)} \in W_1$, $w_{(2)} \in W_2$ and $n \in \C$,
\[
\pi_{n-m}(I(w_{(1)}\otimes w_{(2)})) =0 \;\;\;\mbox{for}\; m \in \N\;\mbox{sufficiently large},
\]
where $\pi_n$ is the canonical projection of $\overline{W}$ to the weight subspace $W_{(n)}$; the {\em Jacobi identity}:
\begin{eqnarray}
\lefteqn{\dps x_{0}^{-1}\delta \bigg( {x_1-z\over x_{0}}\bigg)
Y_3(v,x_1)I(w_{(1)}\otimes w_{(2)})}\nn
&&= z^{-1}\delta \bigg( {x_1-x_{0}\over z}\bigg)
I(Y_1(v,x_1)w_{(1)}\otimes w_{(2)})\nn
&&\hspace{2em}{\dps +x_0^{-1}\delta \bigg( {z-x_1\over -x_0}\bigg)I(w_{(1)}\otimes Y_2(v,x_1)w_{(2)})}
\end{eqnarray}
for $v\in V$, $w_{(1)}\in W_1$ and $w_{(2)}\in W_2$.}
\end{defn}

\begin{rema}{\rm The vector spaces of $P(z)$-intertwining maps of type ${W_3\choose W_1\,W_2}$ is denoted by $\mathcal{M}[P(z)]_{W_1 W_2}^{W_3}$. There is a linear isomorphism between this spaces and the space of intertwining operators of the same type. In particular,
\[
\mbox{dim}\; \mathcal{M}[P(z)]_{W_1 W_2}^{W_3}= \mathcal{N}_{W_1 W_2}^{W_3}.
\]
}
\end{rema}

\begin{defn}{\rm Let $W_1$ and $W_2$ be $V$-modules. A {\em $P(z)$-product} of $W_1$ and $W_2$ is a $V$-module $(W_3, Y_3)$ together with a $P(z)$-intertwining map $I_3$ of type ${W_3\choose W_1\,W_2}$. We denote it by $(W_3, Y_3; I_3)$ or simply by $(W_3, I_3)$. Let $(W_4, Y_4; I_4)$ be another $P(z)$-product of $W_1$ and $W_2$. A {\em morphism} from $(W_3, Y_3; I_3)$ to $(W_4, Y_4; I_4)$ is a module map $\eta$ from $W_3$ to $W_4$ such that
\[
I_4 = \bar{\eta}\circ I_3,
\]
where $\bar{\eta}$ is the natural map from $\overline{W_3}$ to $\overline{W_4}$ uniquely extending $\eta$.}
\end{defn}

Let $\mathcal{C}$ be the category of $V$-modules. The notion of $P(z)$-tensor products of $W_1$ and $W_2$ in $\mathcal{C}$ are defined in terms of a universal property as follows:
\begin{defn}{\rm For $W_1, W_2 \in \mbox{ob}\;\mathcal{C}$, a {\em $P(z)$-tensor product}  of $W_1$ and $W_2$ in $\mathcal{C}$ is a $P(z)$-product  $(W_0, Y_0; I_0)$ with $W_0 \in \mbox{ob}\;\mathcal{C}$ such that for any $P(z)$-product  $(W, Y; I)$ with $W \in \mbox{ob}\;\mathcal{C}$, there is a unique morphism from $(W_0, Y_0; I_0)$ to $(W, Y; I)$. Clearly, a $P(z)$-tensor product  of $W_1$ and $W_2$ in $\mathcal{C}$, if it exists, is unique up to isomorphisms. We denote the $P(z)$-   $(W_0, Y_0; I_0)$ by
\[
(W_1 \boxtimes_{P(z)} W_2, Y_{P(z)}; \boxtimes_{P(z)})
\]
and call the object
\[
(W_1 \boxtimes_{P(z)} W_2, Y_{P(z)})
\]
the {\em $P(z)$-tensor product module}.
}
\end{defn}

The following theorem gives the existence of $P(z)$-tensor product for finitely reductive vertex operator algebras:
\begin{thm}[\cite{HL3, HLZ3}]
Let $V$ be a finitely reductive vertex operator algebra and let $W_1, W_2$ be two $V$-modules. Then $(W_1 \boxtimes_{P(z)}W_2, Y_{P(z)}; \boxtimes_{P(z)})$ exists, and we have a $V$-module decompositions
$$
W_1 \boxtimes_{P(z)}W_2 \simeq\coprod_{i=1}^k (\mathcal{M}[P(z)]_{W_1 W_2}^{M_i})^*\otimes M_i,$$
where $\{M_1, \dots, M_k\}$ is a set of representatives of the equivalence classes of irreducible $V$-modules.
\end{thm}

From this theorem and Theorem \ref{fin-red}, we obtain immediately:
\begin{cor}
If $\ell$ is admissible, then $(W_1 \boxtimes_{P(z)}W_2, Y_{P(z)}; \boxtimes_{P(z)})$ exists for objects $W_1$ and $W_2$ of $\mathcal{O}_{\ell, \ord}$.
\end{cor}

We now compare the contragredient module of $W_1 \circ_{P(z)} W_2$ and a $P(z)$-tensor product module of $W_1$ and $W_2$.
We first recall the construction of $P(z)$-tensor product in \cite{HLZ4}. Let $v \in V$ and
\[
Y_t(v, x) = \sum_{n \in Z}(v \otimes t^n)x^{-n-1} \in (V \otimes \C[t, t^{-1}])[[x, x^{-1}]].
\]
Denote by $\tau_{P(z)}$ the action of
\[
V \otimes \iota_{+}{\mathbb C}[t,t^{- 1}, (z^{-1}-t)^{-1}]
\]
on the vector space $(W_1 \otimes W_2)^*$ given by
\begin{eqnarray}\label{eq41}
(\tau_{P(z)}(\xi)\lambda)(w_{(1)}\otimes w_{(2)})&=&\lambda
(\tau_{W_{1}}((\iota_+\circ T_z\circ\iota_-^{-1}\circ o)\xi)w_{(1)}
\otimes w_{(2)})\nn\\
&&+\lambda (w_{(1)}\otimes\tau_{W_{2}}((\iota_+\circ
\iota_-^{-1}\circ o)\xi)w_{(2)})
\end{eqnarray}
for $\lambda \in (W_1 \otimes W_2)^*$, $w_{(1)} \in W_1$ and $w_{(2)} \in W_2$. Denote by $Y_{P(z)}'$ the action of $V \otimes \C[t,t^{-1}]$ on $(W_1\otimes W_2)^*$ defined by
\[
Y_{P(z)}'(v,x) = \tau_{P(z)}(Y_t(v,x)).
\]
We also have the operators $L_{P(z)}'(n)$ for $n \in \Z$ defined by
\begin{equation}\label{eq43}
Y_{P(z)}'(\omega,x) = \sum_{n \in \Z}L_{P(z)}'(n)x^{-n-2}.
\end{equation}

Given two $V$-modules $W_1$ and $W_2$, let $W_1\hboxtr_{P(z)}W_2$ be the vector space consisting of all the elements $\lambda \in (W_1 \otimes W_2)^*$ satisfying the following two conditions:
\begin{itemize}
\item[(1)]{\it $P(z)$-compatibility condition}:
\begin{itemize}
\item[(a)]{\it Lower truncation condition}: For all $v \in V$, the formal Laurent series $Y_{Q(z)}'(v,x)\lambda$ involves only finitely many negative powers of $x$.
\item[(b)] The following formula holds:
\begin{eqnarray}\label{eq44}
&&\tau_{P(z)}\bigg(z^{-1}\delta\bigg(\frac{x_1-x_0}{z}\bigg)Y_t(v,x_0)\bigg)\lambda\nn
&& \;\;\; = z^{-1}\delta\bigg(\frac{x_1-x_0}{z}\bigg)Y_{P(z)}'(v,x_0)\lambda \;\;\;\mbox{for all} \; v \in V.
\end{eqnarray}
\end{itemize}
\item[(2)] {\it $P(z)$-local grading restriction condition}:
\begin{itemize}
\item[(a)] {\it Grading condition}: $\lambda$ is a (finite) sum of generalized eigenvectors of $(W_1 \otimes W_2)^*$ for the operator $L_{Q(z)}'(0)$.
\item[(b)]The smallest subspace $W_{\lambda}$ of $(W_1 \otimes W_2)^*$ containing $\lambda$ and stable under the component operators $\tau_{P(z)}(v\otimes t^n)$ of the operators $Y_{P(z)}'(v,x)$ for $v \in V$, $n \in \Z$, have the properties:
    \begin{equation*}
    \mbox{dim}\; (W_{\lambda})_{[n]} < \infty
    \end{equation*}
    \[
    (W_{\lambda})_{[n+k]} = 0\;\;\;\mbox{for $k \in \Z$ sufficiently negative};
    \]
    for any $n \in \C$, where the subscripts denote the $\C$-grading by $L_{P(z)}'(0)$-eigenvalues.
\end{itemize}
\end{itemize}

The following theorem gives the construction of $P(z)$-tensor product:

\begin{thm}[\cite{HLZ4}]
The vector space $W_1\hboxtr_{P(z)}W_2$ is closed under the action $Y_{P(z)}'$ of $V$ and the Jacobi identity holds on $W_1\hboxtr_{P(z)}W_2$. Furthermore, the $P(z)$-tensor product of $W_1, W_2 \in \mathcal{C}$ exists if and only if $W_1\hboxtr_{P(z)}W_2$ equipped with $Y_{P(z)}'$ is an object of $\mathcal{C}$ and in this case, the $P(z)$-tensor product is the contragredient module of $(W_1\hboxtr_{P(z)}W_2, Y_{P(z)}')$.
\end{thm}

In the case that $\ell \notin \Q_{\geq 0}$, Zhang proved in \cite{Zh} that for two modules $W_1$ and $W_2$ in the category studied in \cite{KL2}, the vector space $W_1\circ_{Q(z)}W_2$ constructed in \cite{KL2} and the $Q(z)$-tensor product $W_1\hboxtr_{Q(z)}W_2$ for $Q(z)$, the Riemann sphere  equipped with different punctures and local coordinates from those for $P(z)$, are the same. The proof does not use the hypothesis on $\ell$, and in particular, works for the case when $\ell$ is admissible. The proof in fact also works for $P(z)$. Thus we have:
\begin{thm}
For $W_1$ and $W_2$ in $\mathcal{O}_{\ell, \ord}$, the two subspaces $W_1 \circ_{P(z)}W_2$ and $W_1\hboxtr_{P(z)}W_2$ are equal to each other. In particular, $W_1\hboxtr_{P(z)}W_2$ equipped with $Y_{P(z)}'$ is an object of $\mathcal{O}_{\ell, \ord}$, and the $P(z)$-tensor product exists in $\mathcal{O}_{\ell, \ord}$.
\end{thm}

We omit the proof of this result since it is not needed in this paper. As is mentioned above, the proof is similar to the one in \cite{Zh} for $Q(z)$.

\setcounter{equation}{0}
\section{Braided tensor category structures on category $\mathcal{O}_{\ell, \ord}$}

In this section, we verify that the convergence condition and the  expansion condition for intertwining maps (see \cite{HLZ5} and \cite{HLZ6}) hold and 
therefore by the main theorem in \cite{HLZ8}, the category $\mathcal{O}_{\ell, \ord}$ has a natural braided tensor category structure with a twist.
We achieve this by  showing that the intertwining operators among the modules in $\mathcal{O}_{\ell, \ord}$ for an admissible 
number $\ell$ satisfying the assumptions in \cite{HLZ7} and \cite{H4} needed for them to have the associativity property.

The following convergence and extension property was introduced in \cite{HLZ7}:

\begin{defn}
{\rm Let $V$ be a vertex operator algebra. We say that the product of intertwining operators $\mathcal{Y}_1$ and $\mathcal{Y}_2$ satisfies the {\em convergence and extension property for products} if for any $\beta \in \tilde{A}$, there exists an integer $N_{\beta}$ depending only on $\mathcal{Y}_1$ and $\mathcal{Y}_2$ and $\beta$, and for any doubly homogeneous elements $w_{(1)} \in W_1^{(\beta_1)}$ and $w_{(2)} \in W_2^{(\beta_2)}$ $(\beta_1, \beta_2 \in \tilde{A})$ and any $w_{(3)} \in W_3$ and $w_{(4)}' \in W_4'$ such that \[\beta_1 + \beta_2 = -\beta,\]there exist $M \in \N$, $r_k, s_k \in \R$, $i_k, j_k \in \N$, $k = 1, \dots, M$, and analytic functions $f_k(z)$ on $|z| < 1$, $k = 1, \dots, M$, satisfying
\begin{equation}\label{eq70}
\wt w_{(1)} + \wt w_{(2)} + s_k > N_{\beta},\ k = 1, \dots, M,
\end{equation}
such that
\begin{equation}\label{eq71}
\langle w_{(4)}', \mathcal{Y}_1(w_{(1)},  x_1)\mathcal{Y}_2(w_{(2)}, x_2)w_{(3)} \rangle_{W_4} |_{x_1 = z_1,\ x_2 = z_2}
\end{equation}
is absolutely convergent when $|z_1| > |z_2| > 0$ and can be analytically extended to the multivalued analytic function
\begin{equation}\label{eq72}
\sum_{k = 1}^M z_2^{r_k}(z_1 - z_2)^{s_k}(\log z_2)^{i_k}(\log (z_1 - z_2))^{j_k}f_k(\frac{z_1 - z_2}{z_2})
\end{equation}
in the region $|z_2| > |z_1 - z_2| > 0$.

When $i_k = j_k = 0$ for $k = 1, \dots, M$, we call the property above the {\it convergence and extension property without logarithms for products}.
}\end{defn}

The convergence and extension property is part of the sufficient condition for the existence of associativity isomorphism:
\begin{thm}[Theorem 10.3 \cite{HLZ6}, Theorem 11.4 \cite{HLZ7}]\label{assiso}
Assume that $V$ satisfies the following conditions:
\begin{enumerate}
\item[(1)]Every finitely-generated lower-bounded doubly graded generalized $V$-module is a $V$-module.
\item[(2)]The product or the iterates of the intertwining operators for $V$ have the convergence and extension property.
\end{enumerate}
Then for any $V$-module $W_1, W_2$ and $W_3$ and any complex number $z_1$ and $z_2$ satisfying $|z_1|>|z_2|>|z_1-z_2|>0$, there is a unique isomorphism $\mathcal{A}_{P(z_1), P(z_2)}^{P(z_1-z_2), P(z_2)}$ from $W_1 \boxtimes_{P(z_1)} (W_2 \boxtimes_{P(z_2)} W_3)$ to $(W_1 \boxtimes_{P(z_1-z_2)} W_2) \boxtimes_{P(z_2)} W_3$ such that for $w_{(1)} \in W_1$, $w_{(2)} \in W_2$, $w_{(3)} \in W_3$,
\begin{eqnarray*}
&&\overline{\mathcal{A}}_{P(z_1), P(z_2)}^{P(z_1-z_2), P(z_2)}(w_1 \boxtimes_{P(z_1)} (w_2 \boxtimes_{P(z_2)} w_3))\nn
&& \;\;\; = (w_1 \boxtimes_{P(z_1-z_2)} w_2) \boxtimes_{P(z_2)} w_3,
\end{eqnarray*}
where
\[
\overline{\mathcal{A}}_{P(z_1), P(z_2)}^{P(z_1-z_2), P(z_2)}: \overline{W_1 \boxtimes_{P(z_1)} (W_2 \boxtimes_{P(z_2)} W_3)} \rightarrow \overline{(W_1 \boxtimes_{P(z_1-z_2)} W_2) \boxtimes_{P(z_2)} W_3}
\]
is the canonical extension of $\mathcal{A}_{P(z_1), P(z_2)}^{P(z_1-z_2), P(z_2)}$.
\end{thm}

It was observed in \cite{H4} (see Theorem 3.1 in \cite{H4}) that Condition $(1)$ in Theorem \ref{assiso} can in fact be replaced by the following  much weaker condition:
\begin{itemize}
\item [(1)'] For any two modules $W_1$ and $W_2$ and any $z \in \C^{\times}$, if the generalized $V$-module $W_{\lambda}$ generated by a generalized eigenvector $\lambda \in (W_1 \otimes W_2)^{*}$ for $L_{P(z)}'(0)$ satisfying the $P(z)$-compatibility condition is lower truncated, then $W_{\lambda}$ is a $V$-module.
\end{itemize}

We now prove Condition $(1)'$ for $\mathcal{O}_{\ell,\ord}$:
\begin{thm}
Let $\ell$ be an admissible number, $W_1, W_2 \in \mathcal{O}_{\ell, \ord}$, and $\lambda \in (W_1 \otimes W_2)^*$ be a generalized eigenvector for $L_{P(z)}'(0)$ satisfying the $P(z)$-compatibility condition. If the generalized $V$-module $W_{\lambda}$ generated by $\lambda$ is lower truncated, then $W_{\lambda} \in \mathcal{O}_{\ell,\ord}$.
\end{thm}
\pf Since $W_{\lambda}$ is lower truncated, by the characterization (\ref{defnKL}) of $W_1 \circ_{P(z)} W_2$, elements of $W_{\lambda}$ must be in $W_1 \circ_{P(z)} W_2$. By Theorem \ref{closedness}, $W_1 \circ_{P(z)} W_2$ is in $\mathcal{O}_{\ell, \ord}$. Thus $W_{\lambda}$ as a submodule of $W_1 \circ_{P(z)} W_2$ is also in $\mathcal{O}_{\ell, \ord}$. \epfv

The convergence and extension property for $\mathcal{O}_{\ell, \ord}$ is an application of the results in \cite{H3}:
\begin{prop}\label{convext}
If $\ell$ is an admissible number, then products of intertwining operators for $L(\ell,0)$ have the convergence and extension property.
\end{prop}
\pf Since the objects in $\mathcal{O}_{\ell, \ord}$ satisfy the $C_1$-cofiniteness condition, by Theorem 1.4 and Theorem 2.3 in \cite{H3}, the matrix elements of the products and iterates of intertwining operators satisfy a system of differential equations with regular singular points. Thus products and iterates of intertwining operators satisfy the convergence and extension property. \epfv

\begin{rema}
{\rm Since objects in $\mathcal{O}_{\ell, \ord}$ are (ordinary) modules for $L(\ell,0)$, logarithmic intertwining operators are intertwining operators (without logarithmic terms). Thus the associativity of logarithmic intertwining operators in fact gives the associativity of intertwining operators (without logaritmic terms). In partciular, we see that the convergence and extension property (without logarithm) introduced in \cite{H1} holds. } 
\end{rema}

We refer the reader to the books \cite{Tu, BK, EGNO} for basic notions in the theory of tensor categories. 

As a consequence of Theorem \ref{assiso}--\ref{convext} and Corollary 12.16 in \cite{HLZ8}, we obtain the main result of this paper:
\begin{thm}\label{main}
If $\ell$ is an admissible number, then the category $\mathcal{O}_{\ell, \ord}$ has a natural structure of braided tensor category with the tensor product bifunctor $\boxtimes_{P(1)}$, the unit object $L(\ell,0)$ and the braiding isomorphism, the associativity isomorphism, the left and right unit isomorphisms given in Section 12.2 in \cite{HLZ8}. 
\end{thm}

The proof of the following result is the same as the proof of Theorem 4.1 in \cite{H6}:

\begin{prop}\label{twist}
The braided tensor category given in Theorem \ref{main} has a twist given by $e^{2\pi i L_0}$ on objects in $\mathcal{O}_{\ell, \ord}$. 
\end{prop}

\begin{rema}
{\rm As is discussed in the introductory part of Section 12 in \cite{HLZ8}, the results in \cite{HLZ8} in fact give a vertex tensor category sructure in the sense of \cite{HL2} on the category $\mathcal{O}_{\ell, \ord}$. The braided tensor category structure and the twist obtained in Theorem \ref{main} are then induced from the vertex tensor category structure following a standard procedure. }
\end{rema}

\setcounter{equation}{0}
\section{Ribbon categories in the case  $\mathfrak g=\mathfrak{sl}_2$}

As an application of our results above we will now show that in the instance of $\mathfrak{sl}_2$ we get a fusion category and sometimes even a modular tensor category. Recall that a fusion category is a semisimple braided tensor category with only finitely many inequivalent simple objects, such that every object is rigid. A ribbon category is a fusion category with a ribbon twist (which in the instance of a vertex algebra is just given by the action of $e^{2\pi i L_0}$ in Proposition \ref{twist}). Let $(\mathcal C, \boxtimes)$ be a ribbon category over $\mathbb C$ and $X, Y$ be objects in $\mathcal C$, then rigidity implies that there exists a trace on objects, $\text{tr}_X: \text{End}(X) \rightarrow \mathbb C$, while braiding is a natural family of commutativity isomorphisms $c_{X, Y}: X \boxtimes Y \rightarrow Y \boxtimes X$.
Finally, a ribbon category is modular if the matrix of the Hopf link invariants
\[
S_{X, Y} = \text{tr}_{X\boxtimes Y} \left( c_{Y, X} \circ c_{X, Y} \right)
\]
  is non-degenerate. The dimension of an object $X$ is of course the trace of the identity on $X$. We refer the reader to the books \cite{Tu, BK, EGNO} for details on tensor categories. 
 
Let  $\mathfrak g=\mathfrak{sl}_2$ and let $\ell = -2 + \frac{a}{b}$ for $a, b$ coprime positive integers and $a \neq 1$. We will prove that $\mathcal{O}_{\ell, \ord}$ is always ribbon and modular if and only if $b$ is odd. 
The idea is to prove that this tensor category is  braided equivalent to a full fusion subcategory of a rational Virasoro vertex operator algebra. 

Let us start by recalling a few properties of the rational Virasoro vertex operator algebra. The standard reference is the book \cite{IK}. Let $a, b\in \Z_{>1}=\Z_{+}+1$ be coprime and denote by $\Vir{a, b}$ the simple Virasoro vertex algebra at central charge $1-6\frac{(a-b)^2}{ab}$.  Then $\Vir{a, b}$ is $C_2$-cofinite and 
 every $\N$-gradable weak $\Vir{a, b}$-module is a direct sum of irreducible $\Vir{a, b}$-modules (see \cite{Wang}). 
Simple modules $M_{r, s}$ are parameterized by integers $r, s$ in the range $1\leq r \leq a-1, 1\leq s \leq b-1$ and subject to the identification $M_{r, s} \cong M_{a-r, b-s}$. It is clear from the fusion rules (see below)  that $M_{1,1}$ is the identity and each $M_{r, s}$ is its own dual under the fusion product $\boxtimes_{P(1)}$. Then the category of  $\Vir{a, b}$-modules has a natural structure of a modular tensor category (see \cite{Wang}, \cite{H6})
and we denote it by $\mathcal C(\Vir{a, b})$. 

For the fusion rules, let $c$ be a positive integer and $t, t', t''$ integers. Then the numbers
\[
N^c_{t, t', t''} := \begin{cases} 1  & \text{if}\ |t-t'|+1 \leq t'' \leq \text{min}\{t+t'-1, 2c-t-t'-1\}\  \text{and}\\ &  t+t'+t'' \ \text{odd}     \\ 0 & \text{else} \end{cases}
\]
specify the fusion rules, i.e. 
\begin{equation}\label{fl}
M_{r,s } \boxtimes_{P(1)}M_{r', s'} \cong  \coprod_{r''=1}^{a-1} \coprod_{s''=1}^{b-1}  N^a_{r, r', r''}N^b_{s, s', s''} M_{r'', s''}.
\end{equation}
We also need the normalized Hopf link invariants
\[
\frac{S_{(r, s)(r',s')}}{S_{(1, 1)(r', s')}} := \frac{\text{tr}_{M_{r,s } \boxtimes_{P(1)}M_{r', s'}}\left(  c_{M_{r', s'}, M_{r, s}}\circ c_{M_{r, s}, M_{r', s'}} \right)}{\text{dim}\left(M_{r', s'}\right)} ,
\]
where $c_{M_{r, s}, M_{r', s'}}: M_{r,s } \boxtimes_{P(1)}M_{r', s'} \rightarrow M_{r', s' } \boxtimes_{P(1)}M_{r, s}$ is the braiding isomorphism. 
By Verlinde's formula \cite{H5, H6} this coincides with the normalized modular $S$-matrix which gives
\[
\frac{S_{(r, s)(r',s')}}{S_{(1, 1)(r', s')}} = (-1)^{(r+1)s'+(s+1)r'} \frac{\sin\left(\frac{\pi b}{a}rr' \right)\sin\left(\frac{\pi a}{b}ss' \right)}{\sin\left(\frac{\pi b}{a}r' \right)\sin\left(\frac{\pi a}{b}s' \right)}.
\]
Now let $\mathcal C_a$ be the full tensor subcategory whose simple modules are $M_{r, 1}$ with $1\leq r \leq a-1$. Then $\mathcal C_a$ is a ribbon category as each $M_{r, 1}$ is its own dual. We now wonder if it can even be a modular tensor category. 
Brugui\`eres \cite{Bru} has found good criteria to answer this question. Namely, 
\begin{criterion}\label{crit:bru}
A ribbon category $\mathcal C$ is modular if and only if the tensor identity is the only simple object $X$ that satisfies
\begin{equation}\label{eq:forbiddenrelation}
S_{X, Y} = \text{dim}(X)\text{dim}(Y) \qquad \forall \ Y.
\end{equation} 
\end{criterion}
We are thus looking for simple modules $X$ such that
\[
\frac{S_{X, Y}}{\dim(Y)} 
\]
is independent of $Y$. If we cannot find such $X$ that is also not the tensor identity, then $\mathcal C$ is modular.
Moreover, Brugui\`eres \cite{Bru}  essentially proved the following:
\begin{prop}\label{prop:scalarcriterion}
Let $\mathcal C$ be a ribbon category, s.t. for any simple object $X$ the ribbon twist $\theta_X=\widetilde\theta_X \mbox{Id}_X$ with $|\widetilde \theta_X|=1$. Let $X, Y$ be simple objects of $\mathcal C$, then $\theta_{X \boxtimes Y}$ is a scalar if and only if $|S_{X, Y}| = \dim(X)\dim(Y)$.
\end{prop}
\pf
Recall the balancing axiom of a ribbon category
\[
\theta_{X \boxtimes Y} = c_{Y, X} \circ c_{X, Y} \circ \left(\theta_X \boxtimes \theta_Y\right).
\]
Assuming that $\theta_{X\boxtimes Y}$ is a scalar and  then taking the absolute value of the trace of both sides, one gets the identity
\[
| \text{dim}(X\boxtimes Y)| = |S_{X, Y}|.
\]
But $\dim(X\boxtimes Y) = \dim(X)\dim(Y)$.
Conversely, assume that $|S_{X, Y}| = \dim(X)\dim(Y)$. The balancing axiom implies that $|S_{X, Y}|=|\tr_{X\boxtimes Y}\left(\theta_{X\boxtimes Y}\right)|$. So that by Cauchy-Schwarz inequality this is only possible if the twist acts by the same scalar on each simple summand of $X\boxtimes Y$. 
\epfv

Consider $M_{r, 1}$ with $r\neq 1, a-1$, then from (\ref{fl}) the fusion product with the fundamental representation is
\[
M_{2, 1} \boxtimes_{P(1)} M_{r, 1} \cong M_{r-1, 1} \oplus M_{r+1, 1}.
\]
Now, the ribbon twist is given by $e^{2\pi i L_0}$ on a simple module so that we only need to know the conformal weight $h$ of the top level space of the modules. For $M_{r, 1}$ it is
\[
h_{r, 1} = \frac{(rb-a)^2-(a-b)^2}{4ab}
\]
so that we see that
\[
h_{r-1, 1}-h_{r+1, 1} = \frac{(b(r-1)-a)^2-(b(r+1)-a)^2}{4ab} = -\frac{rb}{a} \mod \Z.
\]
It follows that besides the tensor identity $M_{1, 1}$ only the simple current $X=M_{a-1, 1}$ can possibly satisfy \eqref{eq:forbiddenrelation}. 
We thus verify
\begin{equation}\label{eq:hopflinksimplecurrent}
\begin{split}
\frac{S_{(a-1, 1)(r',1)}}{S_{(1, 1)(r', 1)}} &= (-1)^{a} \frac{\sin\left(\frac{\pi b}{a}(a-1)r' \right)\sin\left(\frac{\pi a}{b} \right)}{\sin\left(\frac{\pi b}{a}r' \right)\sin\left(\frac{\pi a}{b} \right)} 
= (-1)^{a+br'+1}
\end{split}
\end{equation}
It follows that 
\[
\frac{S_{(a-1, 1)(r',1)}}{S_{(1, 1)(r', 1)}} 
\]
is independent of $r'$ if and only if $b$ is even. Recall that the categorical dimension is the Hopf link invariant $\text{dim}(X) = S_{X, {\bf 1}}$ with ${\bf 1}$ the tensor identity. 
Thus
\[
S_{(a-1, 1)(r',1)} = \dim (M_{a-1,1})\dim (M_{r',1})
\]
for all $1 \leq r' \leq a-1$ if and only if $b$ is even.
By Criterion \ref{crit:bru}, we obtain:
\begin{thm}\label{thm:virmtc}
The subcategory $\mathcal C_a$ of $\mathcal C(\Vir{a, b})$ is a modular tensor category if and only if $b$ is odd.
\end{thm}
We now use the theory of vertex operator algebra extensions \cite{KO, HKL, CKM} to prove an equivalence of $\mathcal C_a$ to a category closely related to the category of ordinary modules $\mathcal O_{\ell, \ord}$ of $L(\ell, 0)$ and $\mathfrak g=\mathfrak{sl}_2$. Let $\ell=-2+\frac{a}{b}$ with $a, b$ coprime positive integers and $a>1$. Consider the subcategory $\widetilde{\mathcal O}_{\ell, \ord}$ of  $\mathcal O_{\ell, \ord} \otimes \mathcal O_{1, \ord}$ whose simple objects are $L(\ell, \lambda)\otimes L(1, \nu)$ such that $\lambda= \nu \mod A_1$. Here we denote the highest weight module of a highest weight $\lambda$ at level $\ell$  by $L(\ell, \lambda)$ and the root lattice of $\mathfrak{sl}_2$ by $A_1$. 
\begin{thm}\label{tilde}
The categories $\widetilde{\mathcal{O}}_{\ell, \ord}$  and $\mathcal{C}_a$
are equivalent as braided tensor categories with twists. In particular, 
$\widetilde{\mathcal{O}}_{\ell, \ord}$ is rigid and thus is a ribbon category.
\end{thm}
\pf
It is well-known (see Remark 10.3 of \cite{IK}  and \cite{ACL} for a recent proof) that $\text{Com}\left(L(\ell+1, 0), L(\ell, 0)\otimes L(1, 0)\right) \cong \Vir{a, a+b}$. Let $\omega$ be the fundamental weight of $\mathfrak{sl}_2$. 
Then the branching rules are
\begin{equation}
L(\ell, r\omega) \otimes L(1, t\omega) \cong \bigoplus_{\substack{s = 0\\ s + t+r \ \text{even}}}^{a+b-2} L(\ell+1, s\omega) \otimes M_{r+1, s+1}
\end{equation}
for $r=0, \dots, a-2$ and $t=0, 1$. 
We are thus in the situation of \cite{CKM} that $A:=L(\ell, 0)\otimes L(1, 0)$ is a vertex operator algebra extension of $L(\ell+1, 0) \otimes \Vir{a, a+b}$ and there is an induced functor  
$\mathcal{F}$
that relates modules of the two corresponding vertex tensor categories. 
It is given by 
\[
\mathcal F: \mathcal O_{\ell+1, \ord} \otimes \mathcal C(\Vir{a, a+b}) \rightarrow \mathcal C_A, \qquad X\mapsto A\boxtimes_{P(1)}X,
\]
where $\mathcal C_A$ is the cateogry of modules for the algebra object $A$ and $\boxtimes_{P(1)}$ is the $P(1)$-tensor product bifunctor in 
the category $\mathcal O_{\ell+1, \ord} \otimes \mathcal C(\Vir{a, a+b})$. By \cite{HKL} objects of $\mathcal C_A$ are modules for the vertex operator algebra $A=L(\ell, 0)\otimes L(1, 0)$ if and only if they are local $A$-modules. Being local means having trivial monodromy with $A$ or equivalently that the twist is acting as an $A$-morphism. 
We consider the full ribbon subcategory of $\mathcal O_{\ell+1, \ord} \otimes \mathcal C(\Vir{a, a+b})$ whose simple objects are $L(\ell+1, 0) \otimes M_{r, 1}$. This category is clearly ribbon equivalent to $\mathcal C_a$. 
Using the Virasoro fusion rules (\ref{fl}) we see that the induced functor satisfies
\begin{equation}\nonumber
\begin{split}
\mathcal F\left(L(\ell+1, 0) \otimes M_{r, 1}\right) & \cong  \bigoplus_{\substack{s = 0\\ s  \ \text{even}}}^{a+b-2} \left(L(\ell+1, s\omega) \otimes M_{1,  s+1}\right) \boxtimes_{P(1)} \left(L(\ell+1, 0) \otimes M_{r, 1}\right) \\
 &\cong \bigoplus_{\substack{s=0\\ s \ \text{even}}}^{a+b-2} L(\ell+1, s\omega) \otimes M_{r, s+1}\\ &\cong L(\ell, (r-1)\omega) \otimes L(1, \overline{(r-1)}\omega)
\end{split}
\end{equation}
with $\overline{r}=0$ if $r$ even and $\overline{r}=1$ otherwise. Hence we have an essentially surjective tensor functor from $\mathcal C_a$ to $\widetilde{\mathcal O}_{\ell, \ord}$.  The twist clearly acts as a scalar on $\mathcal F\left(L(\ell+1, 0) \otimes M_{r, 1}\right)$ so that this is a local $A$-module. But $\mathcal F$ restricted to modules that induce to local $A$-modules is a braided tensor functor by Theorem 2.67 of \cite{CKM}. It is fully faithful by Proposition 8.2 of \cite{CG} (whose proof follows the argument of the proof of Theorem 5.1 of \cite{OS}). Further by Lemma 1.16 of \cite{KO} (see also Proposition 2.77 of \cite{CKM}) the induced functor maps duals to duals and so especially the rigidity of the object $L(\ell+1, 0) \otimes M_{r, 1}$  gives the rigidity of the object $\mathcal F\left(L(\ell+1, 0) \otimes M_{r, 1}\right)$. 
\epfv

Theorem \ref{tilde} has quite a few corollaries. First we show that  for each admissible level $\ell$, the braided tensor category $\mathcal{O}_{\ell, \ord}$ is a ribbon category:
\begin{cor}
The category $\mathcal O_{\ell,\ord}$ of ordinary $L(\ell, 0)$-modules for the affine vertex algebra of $\mathfrak{sl}_2$ at level $\ell=-2+\frac{a}{b}$ with $a, b$ coprime positive integers and $a>1$ is a ribbon category. 
\end{cor}
\pf The category $\mathcal O_{\ell,\ord}$ is a subcategory of $\mathcal O_{\ell,\ord}\otimes \mathcal O_{1, \ord}$ under the identification of $L(\ell, \lambda)$ with $L(\ell, \lambda) \otimes L(1, 0)$. Thus to prove $L(\ell, \lambda)$ is rigid in $\mathcal{O}_{\ell, \ord}$, it suffices to show $L(\ell, \lambda) \otimes L(1,0)$ is rigid in $\mathcal O_{\ell,\ord}\otimes \mathcal O_{1, \ord}$. If $\lambda = 2k \omega$ for some $k \in \N$, then from Theorem \ref{tilde}, $L(\ell, \lambda) \otimes L(1,0)$ is rigid. Otherwise, using the fusion product
\[
L(\ell, \lambda) \otimes L(1, 0) \cong \left(L(\ell, 0) \otimes L(1, \omega) \right) \boxtimes_{P(1)} \left(L(\ell, \lambda) \otimes L(1, \omega) \right).
\]
From Theorem \ref{tilde}, $L(\ell, \lambda) \otimes L(1, \omega)$ is rigid, also since the simple current $L(\ell, 0) \otimes L(1, \omega)$ is rigid, the fusion product  $L(\ell, \lambda) \otimes L(1, 0)$ is also rigid.
\epfv

Secondly, since the induced functor is a tensor functor \cite[Theorem 2.67]{CKM} and since $L(1, \omega) \boxtimes_{P(1)} L(1, \omega) \cong L(1, 0)$, we get the fusion rules in $\mathcal O_{\ell,\ord}$ from those of $\mathcal C_a$:
\begin{cor}
The fusion rules of simple objects  in $\mathcal O_{\ell,\ord}$ are 
\[
L(\ell, (r-1)\omega) \boxtimes_{P(1)} L(\ell, (r'-1)\omega) \cong \bigoplus_{r''=1}^{a-1} N^a_{r, r', r''}L(\ell, (r''-1)\omega).
\]
\end{cor}
Based on modular properties of characters viewed as formal power series of admissible level $L(\ell, 0)$-modules, a Verlinde formula had been conjectured in \cite{CR1, CR2}. The conjectural fusion rules for ordinary modules are then given in Proposition 15 of \cite{CR2}.
\begin{cor}
Let $\ell=-2+\frac{a}{b}$ with $a, b$ coprime positive integers and $a>1$.
The Verlinde conjecture of \cite{CR1, CR2} for fusion rules of the affine vertex algebra of $\mathfrak{sl}_2$ at level $\ell$ is
 true when restricted to the category $\mathcal O_{\ell, \ord}$ of ordinary modules.
\end{cor}
Finally, we have 
\begin{cor}
Let $\ell=-2+\frac{a}{b}$ with $a, b$ coprime positive integers and $a>1$.
The category $\mathcal O_{\ell, \ord}$ of ordinary $L(\ell, 0)$-modules for the affine vertex algebra of $\mathfrak{sl}_2$ at level $\ell$ is a modular tensor category if and only if $b$ is odd.  
\end{cor}
\pf
The reason for this Corollary to be true is essentially the same as Theorem \ref{thm:virmtc}. Firstly, the ribbon twist on $L(\ell, r\omega)$ and $M_{r+1,1}$ coincide if $r$ is even and differ by the twist on $L(1, \omega)$ which is equal to $e^{\pi i/2}$ if $r$ is odd. It follows using Proposition \ref{prop:scalarcriterion} that the only simple module besides $L(\ell, 0)$ that can possibly satisfy Criterion \ref{crit:bru} is the simple current $L(\ell, a\omega)$.

By Theorem 2.89 of \cite{CKM} the normalized Hopf link invariants of the Virasoro modules $(M_{r+1, 1}, M_{r'+1, 1})$ coincide with the normalized Hopf link invariants of the induced  modules 
\begin{equation}\nonumber
\begin{split}
&\left(\mathcal F\left(L(\ell+1, 0) \otimes M_{1, r+1} \right),  \mathcal F\left(L(\ell+1, 0) \otimes M_{1, r'+1} \right)\right) = \\
&\qquad\qquad\qquad\qquad\qquad\qquad\quad = \left(L(\ell, r\omega) \otimes L(1, \overline{r}\omega), L(\ell, r'\omega) \otimes L(1, \overline{r}'\omega)\right)
\end{split}
\end{equation}
 of the extended vertex operator algebra. 
Let us denote the Hopf link invariants of $L(\ell, 0)$-modules $L(\ell, \lambda)$ and $L(\ell, \mu)$ by $S^{(\ell)}_{\lambda, \mu}$.
The Hopf link invariants of $L(1, 0)$-modules are $\frac{S^{(1)}_{0, \omega}}{S^{(1)}_{0, 0}}=1$ and $\frac{S^{(1)}_{\omega, \omega}}{S^{(1)}_{0, \omega}}=-1$ so that it follows from \eqref{eq:hopflinksimplecurrent} that
\[
\frac{S^{(\ell)}_{a\omega, r\omega}}{S^{(\ell)}_{0, r\omega}} (-1)^{ar} = (-1)^{a+(a+b)(r+1)+1}, 
\]
and hence 
\[
\frac{S^{(\ell)}_{a\omega, r\omega}}{S^{(\ell)}_{0, r\omega}}  =(-1)^{(r+1)b+1}.
\]
The Corollary thus follows from Brugui\`eres' Criterion \ref{crit:bru}.
\epfv

\def\refname{\hfil{REFERENCES}}

\noindent {\small \sc Department of Mathematical and Statistical Sciences, University of Alberta,
Edmonton, AB T6G 2G1}

\noindent {\em E-mail address}: creutzig@ualberta.ca

\vspace{2em}

\noindent {\small \sc Department of Mathematics, Rutgers University,
110 Frelinghuysen Rd., Piscataway, NJ 08854-8019}

\noindent {\em E-mail address}: yzhuang@math.rutgers.edu

\vspace{2em}

\noindent{\small \sc Department of Mathematics, Yale University,
10 Hillhouse Ave, New Haven, CT 06511}

\noindent {\em E-mail addess}: jinwei.yang@yale.edu.

\begin{thebibliography}{EGNO}

\bibitem[ACL]{ACL} T. Arakawa, T. Creutzig and A. Linshaw, Coset construction of principal $W$-algebras, in preparation.

\bibitem[ACR]{ACR} J. Auger, T. Creutzig and D. Ridout, Modularity of logarithmic parafermion vertex algebras, arXiv:1704.05168.


\bibitem[Ad1]{Ada} D. Adamovi$\acute{\mbox{c}}$, Some rational vertex algebras, {\em Glas. Mat. Ser. III}, {\bf 29} (1994), 25--40.

\bibitem[Ad2]{Ad2}
  D. Adamovi$\acute{\mbox{c}}$,  A Construction of admissible $A^{(1)}_1$ modules of level -$4/3$, J. Pure Appl. Algebra, 196:119-134, 2005.

\bibitem[AK]{AK} T. Arakawa and K. Kawasetsu, Quasi-lisse vertex algebras and modular linear differential equations, arXiv:1610.05865.

\bibitem[AL]{AL}J. Axtell and K.-H. Lee, Vertex operator algebras associated to type $G$ affine Lie algebras, {\em J. Alg.}, {\bf 337} (2011), 195--223.

\bibitem[AM]{AM}D. Adamovi$\acute{\mbox{c}}$ and A. Milas, Vertex operator algebra associated to modular invariant representations for $A_1^{(1)}$, {\em Math. Res. Lett.}, {\bf 2} (1995), 563--575.

\bibitem[Ar]{A} T. Arakawa, Rationality of admissible affine vertex algebras in the category $\mathcal{O}$, {\em Duke Math. J.}, {\bf 165} (2016), 67--93.

\bibitem[AY]{AY}H. Awata, Y. Yamada, Fusion rules for the fractional level $\widehat{\mathfrak{sl}}(2, \C)$ algebra, {\em Mod. Phys. Lett.} {\bf A7}, 1185 (1992).

\bibitem[B]{Bru} Alain Brugui\`eres. Cat\'egories pr\'emodulaires, modularisations et invariants des vari\'et\'es de
dimension 3. Math. Ann., 316(2):215-236, 2000.

\bibitem[B--vR]{BLL+} Christopher Beem, Madalena Lemos, Pedro Liendo, Wolfger Peelaers, Leonardo Rastelli, and Balt C. van Rees. Infinite chiral symmetry in four dimensions. Comm. Math. Phys., 336(3):1359-1433, 2015.

\bibitem[BF]{BF}D. Bernard, G. Felder, Fock representations and BRST cohomology in $SL(2)$ current
algebra, {\em Commun. Math. Phys.} {\bf 127} (1990), 145--168.

\bibitem[BFM]{BFM}A. Beilinson, B. Feigin and B. Mazur, Introduction to algebraic field theory on curves, preprint, 1991.

\bibitem[BK]{BK}
B. Bakalov and A. Kirillov, Jr., {\em Lectures
on tensor categories and modular functors},
University Lecture Series, Vol. 21,
Amer. Math. Soc., Providence, RI, 2001.

\bibitem[BPZ]{BPZ}A. Belavin, A. Polyakov, A. Zamolodchikov, Infinite conformal symmetry in two-dimensional quantum field theory, {\em Nuclear Phys. B} {\bf 241} (1984), 333--380.

\bibitem[BR]{BR}  C. Beem and L. Rastelli, Vertex operator algebras, Higgs branches, and modular differential equations, arXiv:1707.07679.

\bibitem[CG]{CG} T. Creutzig, D. Gaiotto, Vertex algebras for S-duality, arXiv:1708.00875.

\bibitem[CKL]{CKL} T. Creutzig, S. Kanade and A. Linshaw, Simple current extensions beyond semi-simplicity, arXiv:1511.08754.

\bibitem[CKLR]{CKLR} T. Creutzig, S. Kanade, A. R. Linshaw and D. Ridout, Schur-Weyl Duality for Heisenberg Cosets, arXiv:1611.00305.

\bibitem[CKM]{CKM} T. Creutzig, R. McRae and S. Kanade, Tensor categories for vertex operator superalgebra extensions,  arXiv:1705.05017.

\bibitem[CMR]{CMR}
  T. Creutzig, A. Milas and M. Rupert, Logarithmic Link Invariants of $\overline{U}_q^H(\mathfrak{sl}_2)$ and Asymptotic Dimensions of Singlet Vertex Algebras, arXiv:1605.05634.

\bibitem[CR1]{CR1} T. Creutzig and D. Ridout, Modular Data and Verlinde Formulae for Fractional Level WZW Models I, {\it Nucl. Phys.} {\bf  B865} (2012),  83--114

\bibitem[CR2]{CR2} T. Creutzig and D. Ridout, Modular Data and Verlinde Formulae for Fractional Level WZW Models II, {\it Nucl. Phys.} {\bf B875} (2013) 423--458.

\bibitem[CRW]{CRW} T. Creutzig, D. Ridout and S. Wood, Coset Constructions of Logarithmic (1, p) Models, Lett.\ Math.\ Phys.\  {\bf 104} (2014) 553.

\bibitem[DLM]{DLM} C. Dong, H. Li, G. Mason, Vertex operator algebras associated to admissible representations
of $\widehat{\mathfrak{sl}(2, \C)}$, {\em Commun. Math. Phys.} {\bf 184} (1997), 65--93.


\bibitem[EGNO]{EGNO} P. Etingof, S. Gelaki, D. Nikshych and V. Ostrik, {\it Tensor Categories}, Mathematical Surveys
and Monographs, Vol. 205, American Mathematical Society, Providence, RI, 2015.

\bibitem[Fa]{Fa}
G. Faltings,
A proof for the Verlinde formula,
{\it J. Alg. Geom.} {\bf 3} (1994),  347--374.

\bibitem[FHL]{FHL}
I.~B. Frenkel, Y.-Z. Huang and J.~Lepowsky, {\em On axiomatic
approaches to vertex operator algebras and modules}, Mem. Amer.
Math. Soc. 104, Amer. Math. Soc., Providence, 1993, no. 494
(preprint, 1989).

\bibitem[Fi1]{F1}M. Finkelberg, Fusion categories, Ph.D. thesis, Harvard University, 1993.

\bibitem[Fi2]{F2}M. Finkelberg, An equivalence of fusion categories, {\em Geom. Funct. Anal.} {\bf 6} (1996), 249--267.

\bibitem[Fi3]{F3}M. Finkelberg, Erratum to: An Equivalence of Fusion Categories, {\em Geom. Funct. Anal.} {\bf 23} (2013), 810–811.

\bibitem[FM]{FM}B. Feigin and F. Malikov, Modular functor and representation theory of $\widehat{\mathfrak{sl}}_2$ at a rational level, In {\em Operads: Proceedings of Renaissance Conferences (Hartford, CT/Luminy,
1995)}, Contemp. Math., Vol 202, (American Mathematical Society, Providence, 1997), pp. 357--405.

\bibitem[FZ]{FZ}
I. B. Frenkel and Y.  Zhu, Vertex operator algebras associated to representations of affine and Virasoro
algebras, {\em Duke Math. J.} {\bf 66} (1992),  123--168.

\bibitem[G]{G}M. Gaberdiel, Fusion rules and logarithmic representations of a WZW model at
fractional level, {\em Nucl. Phys.} {\bf B618} (2001), 407--436.

\bibitem[GK]{GK}M. Gorelik and V. Kac, On complete reducibility for infinite-dimensional Lie algebras, {\em Adv. Math.,} {\bf 262 (2)}, (2011), 1911--1972.

\bibitem[H1]{H}Y.-Z. Huang, {\em Two-dimensional conformal field theory and vertex operator algebra}, Progress in Math., 148, Birkh\"{a}user, Boston, 1997.

\bibitem[H2]{H1}Y.-Z. Huang, A theory of tensor products for module categories for a vertex operator algebra, IV, {\em J. Pure. Appl. Alg.}, {\bf 100} (1995), 173--216.

\bibitem[H3]{H2}Y.-Z. Huang, Virasoro vertex operator algebras, (nonmeromorphic) operator product expansion and the tensor product theory, {\em J. Alg.} {\bf 182} (1996), 201--234.

\bibitem[H4]{H3}
Y.-Z. Huang, Differential equations and intertwining operators, {\em Comm. Contemp. Math}. {\bf 7} (2005),
375--400.

\bibitem[H5]{H5}
Y.-Z. Huang,  Vertex operator algebras and the Verlinde conjecture, {\it Comm. Contemp. Math.} {\bf 10} (2008), 103--154.


\bibitem[H6]{H6}
Y.-Z. Huang, Rigidity and modularity of vertex tensor categories, {\it Comm. Contemp. Math.} {\bf 10} (2008), 871--911.


\bibitem[H7]{H4}Y.-Z. Huang, On the applicability of logarithmic tensor category theory, arXiv: $1702.00133$.

\bibitem[HKL]{HKL}
	Y.-Z. Huang, A. Kirillov and J. Lepowsky, Braided tensor categories and extensions of vertex operator algebras, \textit{Comm. Math. Phys.}  \textbf{337}  (2015),  1143-1159.

\bibitem[HL1]{HL1}Y.-Z. Huang and J. Lepowsky, Toward a theory of tensor products for representations of a vertex operator algebra, in: {\em Proc. 20th Internatinal Conference on Differential Geometric Methods in Theoretical Physics, New York, 1991}, ed. S. Catto and A. Rocha, World Scientific, Singapore, 1992, Vol. 1, 344--354.

\bibitem[HL2]{HL2}Y.-Z. Huang and J. Lepowsky, Tensor products of modules for a vertex operator algebra and vertex tensor categories, in: {\em Lie Theory and Geometry, in honor of Bertram Konstant}, ed. R. Brylinski, J.-L. Brylinski, V. Guillemin, V. Kac, Birkh$\ddot{\mbox{a}}$user, Boston, 1994, 349--383.

\bibitem[HL3]{HL3}Y.-Z. Huang and J. Lepowsky, A theory of tensor products for module categories for a vertex operator algebra, I, {\em Selecta Mathematica, New Series}, {\bf 1} (1995), 699--756.

\bibitem[HL4]{HL4}Y.-Z. Huang and J. Lepowsky, A theory of tensor products for module categories for a vertex operator algebra, II, {\em Selecta Mathematica, New Series}, {\bf 1} (1995), 757--786.

\bibitem[HL5]{HL5}Y.-Z. Huang and J. Lepowsky, A theory of tensor products for module categories for a vertex operator algebra, III, {\em J. Pure. Appl. Alg.}, {\bf 100} (1995), 141--171.

\bibitem[HL6]{HL6}Y.-Z. Huang and J. Lepowsky, A theory of tensor products for module categories for a vertex operator algebra, V, to appear.

\bibitem[HL7]{HL7}Y.-Z. Huang and J. Lepowsky, Intertwining operator algebras and vertex tensor categories for affine Lie algebras, {\em Duke Math. J.} {\bf 99} (1999), 113--134.



\bibitem[HLZ1]{HLZ0}
Y.-Z. Huang, J, Lepowsky and L. Zhang,  A logarithmic generalization of tensor product theory for modules for a vertex operator algebra, 
{\it Internat. J. Math.} {\bf 17} (2006), 975--1012. 

\bibitem[HLZ2]{HLZ1}
Y.-Z. Huang, J, Lepowsky and L. Zhang, Logarithmic tensor category theory for generalized modules for a conformal vertex algebra, I: Introduction and strongly graded algebras and their generalized modules, {\em Conformal Field Theories and Tensor Categories, Proceedings of a Workshop Held at Beijing International Center for Mathematics Research}, ed. C. Bai, J. Fuchs, Y.-Z. Huang, L. Kong, I. Runkel and C. Schweigert, Mathematical Lectures from Beijing University, Vol. 2, Springer, New York, 2014, 169--248.

\bibitem[HLZ3]{HLZ2}
Y.-Z. Huang, J, Lepowsky and L. Zhang, Logarithmic tensor category theory for generalized modules for a conformal vertex algebra, II: Logarithmic formal calculus and properties of logarithmic intertwining operators, arXiv: 1012.4196.

\bibitem[HLZ4]{HLZ3}
Y.-Z. Huang, J, Lepowsky and L. Zhang, Logarithmic tensor category theory for generalized modules for a conformal vertex algebra, III: Intertwining maps and tensor product bifunctors, arXiv: 1012.4197.

\bibitem[HLZ5]{HLZ4}
Y.-Z. Huang, J, Lepowsky and L. Zhang, Logarithmic tensor category theory for generalized modules for a conformal vertex algebra, IV: Construction of tensor product bifunctors and the compatibility conditions, arXiv: 1012.4198.

\bibitem[HLZ6]{HLZ5}
Y.-Z. Huang, J, Lepowsky and L. Zhang, Logarithmic tensor category theory for generalized modules for a conformal vertex algebra, V: Convergence condition for intertwining maps and the corresponding compatibility condition, arXiv: 1012.4199.

\bibitem[HLZ7]{HLZ6}
Y.-Z. Huang, J, Lepowsky and L. Zhang, Logarithmic tensor category theory for generalized modules for a conformal vertex algebra, VI: Expansion condition, associativity of logarithmic intertwining operators, and the associativity isomorphisms, arXiv: 1012.4202.

\bibitem[HLZ8]{HLZ7}
Y.-Z. Huang, J, Lepowsky and L. Zhang, Logarithmic tensor category theory for generalized modules for a conformal vertex algebra, VII: Convergence and extension properties and applications to expansion for intertwining maps, arXiv: 1110.1929.

\bibitem[HLZ9]{HLZ8}
Y.-Z. Huang, J, Lepowsky and L. Zhang, Logarithmic tensor category theory for generalized modules for a conformal vertex algebra, VIII: Braided tensor category structure on categories of generalized modules for a conformal vertex algebra, arXiv: 1110.1931.

\bibitem[HY]{HY}
Y.-Z. Huang and J. Yang, Logarithmic intertwining operators and associative algebras, {\em J. Pure Appl. Alg.}
{\bf 216} (2012), 1467--1492.


\bibitem[IK]{IK} K. Iohara and Y. Koga, {\it Representation Theory of the Virasoro Algebra}, Springer Monographs in Mathematics, Springer-Verlag, London, 2011.

\bibitem[K]{K}
A. W. Knapp, {\em Representation Theory of Semisimple Groups,} Princeton University Press, Princeton, New Jersey, 1986.

\bibitem[KL1]{KL1}D. Kazhdan and G. Lusztig, Affine Lie algebras and quatum groups, {\em International Mathematics Research Notices} (in {\em Duke Math. J.}) {\bf 2} (1991), 21--29.

\bibitem[KL2]{KL2}D. Kazhdan and G. Lusztig, Tensor structure arising from affine Lie algebras, I, {\em J. Amer. Math. Soc.} {\bf 6} (1993), 905--947.

\bibitem[KL3]{KL3}D. Kazhdan and G. Lusztig, Tensor structure arising from affine Lie algebras, II, {\em J. Amer. Math. Soc.} {\bf 6} (1993), 949--1011.

\bibitem[KL4]{KL4}D. Kazhdan and G. Lusztig, Tensor structure arising from affine Lie algebras, III, {\em J. Amer. Math. Soc.} {\bf 7} (1994), 335--381.

\bibitem[KL5]{KL5}D. Kazhdan and G. Lusztig, Tensor structure arising from affine Lie algebras, IV, {\em J. Amer. Math. Soc.} {\bf 7} (1994), 383--453.

\bibitem[KO]{KO} A. Kirillov, Jr., and V. Ostrik, On a $q$-analogue of the McKay correspondence and the $ADE$ classification of $\mathfrak{sl}_2$   conformal field theories, \textit{Adv. Math.}  \textbf{171}  (2002), 183-227. 

\bibitem[KW1]{KW1}Victor G. Kac and M. Wakimoto, Modular invariant representations of infinite-dimensional Lie algebras and superalgebras. {\em Proc. Natl. Acad. Sci. USA} {\bf 85}, 4956--4960 (1988).

\bibitem[KW2]{KW2}Victor G. Kac and M. Wakimoto, Classification of modular invariant representations of affine algebras. In: {\em Infinite-dimensional Lie algebras and groups}, Proceedings of the conference held at CIRM, Luminy, edited by Victor. G. Kac, 1988.

\bibitem[KZ]{KZ}V. Knizhnik and A. Zamolodchikov, Current algebra and Wess-Zumino models in two dimensions, {\em Nuclear Phys. B} {\bf 247} (1984), 83--103.

\bibitem[Li1]{Li1}
H. S. Li, An analogue of the Hom functor and a generalized nuclear democracy
theorem, \textit{Duke Math. J.}  \textbf{93}  (1998),  73--114.

\bibitem[Li2]{Li2}
H. S. Li, Determining fusion rules by $A(V)$-modules and bimodules, \textit{J.
Algebra}  \textbf{212}  (1999),  515--556.



\bibitem[M]{M}A. Milas, Weak modules and logarithmic intertwining operators for vertex operator
algebras, in {\em Recent Developments in Infinite-Dimensional Lie Algebras and Conformal
Field Theory}, eds. S. Berman, P. Fendley, Y.-Z. Huang, K. Misra and B. Parshall,
Contemporary Mathematics, Vol. 297, American Mathematical Society, Providence,
2002, pp. 201--225.

\bibitem[Mc]{Mc}Robert McRae, Non-negative integral level affine Lie algebra tensor categories and their associativity isomorphism, {\em Comm. Math. Phys.} {\bf 346} (2016), 349–395.  

\bibitem[MS]{MS}G. Moore and N. Seiberg, Classical and quatum conformal field theory, {\em Comm. Math. Phys.} {\bf 123} (1989), 177--254.

\bibitem[OS]{OS}V. Ostrik and M. Sun, Level-rank duality via tensor categories, {\it Comm. Math. Phys.} {\bf 326} (2014), 49–61.

\bibitem[P1]{Per1}O. Per$\check{\mbox{s}}$e, Vertex operator algebras associated to type $B$ affine Lie algebras on admissible half-integer levels, {\em J. Algebra}, {\bf 307} (2007), 215--248.

\bibitem[P2]{Per2}O. Per$\check{\mbox{s}}$e, Vertex operator algebras associated to certain admissible modules for affine Lie algebras of type $A$, {\em Glas. Mat. Ser. III}, {\bf 43} (2008), 41--57.

\bibitem[R1]{R1} D. Ridout. Fusion in fractional level sl(2)-theories with $k = -2$ . Nucl. Phys., B848:216-250, 2011.

\bibitem[R2]{R2} D. Ridout. sl(2)$_{-1/2}$ and the Triplet Model,
  Nucl.\ Phys.\ B {\bf 835} (2010) 314.

\bibitem[R3]{R3} D. Ridout. $\widehat{sl}(2)_{-1/2}$: A Case Study,
  Nucl.\ Phys.\ B {\bf 814} (2009) 485.

\bibitem[RW]{RW} D. Ridout and S. Wood, Relaxed singular vectors, Jack symmetric functions and fractional level $\widehat{\mathfrak{sl}}(2)$ models, Nucl.\ Phys.\ B {\bf 894} (2015) 621.

\bibitem[S]{Sawin} S. F. Sawin, Quantum groups at roots of unity and modularity. J. Knot Theory Ramif. 15(10), 1245-1277 (2006).

\bibitem[Te]{T}
C. Teleman, Lie algebra cohomology and the fusion rules,
{\it Comm. Math. Phys.} {\bf 173} (1995), 265--311.

\bibitem[Tu]{Tu}
V. Turaev, {\em Quantum invariants of knots and $3$-manifolds},
de Gruyter Studies in Math., Vol. 18, 
Walter de Gruyter, Berlin, 1994.

\bibitem[TK]{TK}
A. Tsuchiya and Y. Kanie, Vertex operators in conformal field theory on $\mathbf{P}^1$ and monodromy representations of braid group, in \textit{Conformal Field Theory and Solvable Lattice Models (Kyoto, 1986)}, Adv. Stud. Pure Math. \textbf{16}, Academic Press, Boston, 1988, 297-372.

\bibitem[W]{Wang} W. Wang, Rationality of Virasoro vertex operator algebras, {\it Int. Math. Res. Notices} {\bf  1993} (1993), 197-211.

\bibitem[Y1]{Y1}
J. Yang, Some results in the representation theory of strongly graded vertex algebras, Ph.D. thesis, Rutgers
University, May 2014.

\bibitem[Y2]{Y2}
J. Yang, Differential equations and logarithmic intertwining operators for strongly graded vertex algebra,
{\em Commun. Contempt. Math.} {\bf 19} (2017), no.2, 1650009.

\bibitem[Z]{Z}
Y.-C. Zhu, Modular invariance of characters of vertex operator algebras,
\textit{J. Amer. Math. Soc.} \textbf{9} (1996), 237--307.

\bibitem[Zh]{Zh}L. Zhang, Vertex tensor category structure on a category of Kazhdan-Lusztig, {\em New York J. Math.} {\bf 14} (2008), 261--284.
\end{thebibliography}
\end{document}